\documentstyle{amltd}
\begin{document}
\annalsline{157}{2003}
\received{July 25, 1999}
\startingpage{45}
\def\bye{\end{document}}
 \font\tenrm=cmr10

\def\joinrel{\mathrel{\mkern-4mu}}
\def\relbar{\mathrel{\smash-}}
\def\lrar{\relbar\joinrel\relbar\joinrel\relbar\joinrel\relbar\joinrel\rightarrow}

\def\lraop#1{{\displaystyle\mathop{\lrar}^#1}}
\catcode`\@=11
\font\twelvemsb=msbm10 scaled 1100
\font\tenmsb=msbm10
\font\ninemsb=msbm10 scaled 800
\newfam\msbfam
\textfont\msbfam=\twelvemsb  \scriptfont\msbfam=\ninemsb
  \scriptscriptfont\msbfam=\ninemsb
\def\msb@{\hexnumber@\msbfam}
\def\Bbb{\relax\ifmmode\let\next\Bbb@\else
 \def\next{\errmessage{Use \string\Bbb\space only in math
mode}}\fi\next}
\def\Bbb@#1{{\Bbb@@{#1}}}
\def\Bbb@@#1{\fam\msbfam#1}
\catcode`\@=12

 \catcode`\@=11
\font\twelveeuf=eufm10 scaled 1100
\font\teneuf=eufm10
\font\nineeuf=eufm7 scaled 1100
\newfam\euffam
\textfont\euffam=\twelveeuf  \scriptfont\euffam=\teneuf
  \scriptscriptfont\euffam=\nineeuf
\def\euf@{\hexnumber@\euffam}
\def\frak{\relax\ifmmode\let\next\frak@\else
 \def\next{\errmessage{Use \string\frak\space only in math
mode}}\fi\next}
\def\frak@#1{{\frak@@{#1}}}
\def\frak@@#1{\fam\euffam#1}
\catcode`\@=12

 \def\al{\alpha}
\def\b{\betta}
\def\fb{\#_{V}}
\def\fbl{\#_{\la}}
\def\vcup{\cup_{V}}
\def\cu{{\cal C}}
\def\ocu#1{\ov{\cal C}_{#1}}
\def\cn{{\cal N}}
\def\co{{\cal O}}
\def\cj{{\cal J}}
\def\Jn{{\cal J}_\nu}
\def\om{\omega}
\def\bd{\partial}

\def\dim{{\rm dim\; }}
\def\de{\delta}
\def\tiv{\ma\ti_{\ev}}
\def\ga{\gamma}
\def\text#1{{\it #1}}
\def\ca{{\cal A}}
\def\la{\lambda}
\def\sla{\sqrt{\la}}
\def\ep{\varepsilon}
\def\hom{{\rm Hom}}

\def\be{\begin{equation}}
\def\ee{\end{equation}}
\def\bear{\begin{eqnarray}}
\def\eear{\end{eqnarray}}
\def\best{\begin{eqnarray*}}
\def\eest{\end{eqnarray*}}
\def\pf{\demo{Proof}}


\def\sspa{\vphantom{L^N_K}}
\def\mspa{\vphantom{\sum}}
\def\bspa{\vphantom{\int}}

\def\non{\noindent}
\def\pf{\demo{Proof}}

\def\ra{\rightarrow}
\def\lra{\leftarrow}
\def\rg{\rangle}
\def\lg{\langle}
\def\r#1{\right#1}
\def\l#1{\left#1}

\def\ma#1{\mathop {#1} \limits}
\def\fl{{f_\la}}
\def\b{\beta}
\def\bl{\b_\la}

\def\Si{\Sigma}
\def\De{\Delta}

\def\ti{\times}
\def\tit{\ma(\ti)_{T^k}}
\def\tig{\ma(\ti)_G}

\def\del{\overline \partial}


\def\Z{{ \Bbb Z}}
\def\R{{ \Bbb R}}
\def\P{{ \Bbb P}}
\def\Q{{ \Bbb Q}}
\def\cx{{ \Bbb C}}

\def\evl{{{\rm ev}}}
\def\ev{{\rm ev}}

\def\im{{\rm Im}}
\def\wt#1{\widetilde{#1}}
\def\wh#1{\widehat{#1}}
\def\ov#1{\overline{#1}}

\def\w{\omega}
\def\tr{\mbox{tr}}


\def\jh{$(J,\nu)$-holomorphic}
\def\sgn{{\rm sgn}}
\def\tr{{\rm tr}}
\def\M{{\cal M}}
\def\Mdmgk{\M_{g,n}}
\def\bMdmgk{\overline{\M}_{g,n}}
\def\bMdmgkm{\overline{\M}_{g,n,m}}
\def\Aut{{\rm Aut}}
\def\div{{\rm Div}}
\def\J{\cal J}

\def\oC{\ov{\cal C}}
\def\U{{\cal U}}
\def\UM{{\cal UM}}
\def\H{{\cal H}}
\def\ker{{\rm ker }}
\def\cok{{\rm coker }}
\def\ind{{\rm index}}



\title{Relative Gromov-Witten invariants}
\shorttitle{Relative Gromov-Witten invariants}  

 \acknowledgements{The research of both authors was partially supported by
the N.S.F.  The first author was also supported by a Sloan Research
Fellowship.}
 \twoauthors{Eleny-Nicoleta Ionel}{Thomas H. Parker}
 \institutions{University of Wisconsin-Madison,
   Madison, WI\\
{\eightpoint {\it E-mail address\/}: ionel@math.wisc.edu}\\
 \vglue6pt
Michigan
State University, East Lansing,
MI\\
{\eightpoint {\it E-mail address\/}: parker@math.msu.edu
}}


\centerline{\bf Abstract}
\vglue12pt
We define relative Gromov-Witten invariants of a symplectic manifold
relative to a codimension-two symplectic submanifold.  These
invariants are the key ingredients in the symplectic sum formula of
\cite{IP4}. The main step is the construction of a compact space of
`$V$-stable' maps.  Simple special cases include the Hurwitz numbers
for algebraic curves and  the enumerative invariants of Caporaso and
Harris.

\vglue18pt

Gromov-Witten invariants are invariants of a closed symplectic
manifold $(X,\w)$.  To define them, one introduces a compatible almost
complex structure $J$ and a perturbation term $\nu$, and considers the
maps $f:C\to X$ from a genus $g$ complex curve $C$ with $n$ marked
points which satisfy the pseudo-holomorphic map equation $\del f=\nu$
and represent a class $A=[f]\in H_2(X)$. The set of such maps,
together with their limits, forms the compact space of stable maps
$\ov\M_{g,n}(X,A)$.  For each stable map, the domain determines a
point in the Deligne-Mumford space $\ov\M_{g,n}$ of curves, and
evaluation at each marked point determines a point in $X$. Thus there
is a map
\bear
\ov\M_{g,n}(X,A) \to \ov\M_{g,n} \ti X^n.
\label{intro.1}
\eear
The Gromov-Witten invariant of $(X,\w)$ is the homology class of the
image for generic $(J,\nu)$.  It depends only on the isotopy class of
the symplectic structure.  By choosing   bases of the cohomologies of
$\ov\M_{g,n}$ and $X^n$, the GW invariant can be viewed as a
collection of numbers that count the number of stable maps satisfying
constraints.  In important cases these numbers are equal to
enumerative invariants defined by algebraic geometry.

In this article we construct Gromov-Witten invariants for a symplectic
manifold $(X,\w)$ relative to a codimension two symplectic submanifold
$V$.  These invariants are designed for use in formulas describing how
GW invariants behave under symplectic connect sums along $V$ --- an
operation that removes $V$ from $X$ and replaces it with an open
symplectic manifold $Y$ with the symplectic structures matching on the
overlap region.  One expects the stable maps into the sum to be pairs
of stable maps into the two sides which match in the middle.  A sum
formula thus requires a count of stable maps in $X$ that keeps track
of how the curves intersect $V$.

Of course, before speaking of stable maps one must extend $J$ and
$\nu$ to the connect sum.  To ensure that there is such an extension
we require that the pair $(J,\nu)$ be `$V$-compatible' as defined in
Section 3.  For such pairs, $V$ is a\break $J$-holomorphic submanifold ---
something that is not true for generic $(J,\nu)$. The relative
invariant gives counts of stable maps for these special $V$-compatible
pairs.  These counts are different from those associated with the
absolute GW invariants.

The restriction to $V$-compatible $(J,\nu)$ has repercussions. It
means that pseudo-holomorphic maps $f:C\to V$ into $V$ are
automatically pseudo-holo\-morphic maps into $X$.  Thus for
$V$-compatible $(J,\nu)$, stable maps may have domain components whose
image lies entirely in $V$.  This creates problems because such maps
are not transverse to $V$.  Worse, the moduli spaces of such maps can
have dimension {\it larger} than the dimension of $\M_{g,n}(X,A)$.  We
circumvent these difficulties by restricting attention to the stable
maps which {\it have no components mapped entirely into $V$}.  Such
`$V$-regular' maps intersect $V$ in a finite set of points with
multiplicity.  After numbering these points, the space of $V$-regular
maps separates into components labeled by vectors $s=(s_1, \dots ,
s_{\ell})$, where $\ell$ is the number of intersection points and
$s_k$ is the multiplicity of the $k^{\rm th}$ intersection point.  In
Section 4 it is proved that each (irreducible) component
$\M^{V}_{g,n,s}(X,A)$ of $V$-regular stable maps is an orbifold; its
dimension depends on $g,n,A$ and on the vector $s$.

The next step is to construct a space that records the points where a\break
$V$-regular map intersects $V$ and records the homology class of the
map.  There is an obvious map from $\M^V_{g,n,s}(X,A)$ to
$H_2(X) \times V^\ell$ that would seem to serve this purpose.
However, to be
useful for a connect sum gluing theorem, the relative invariant should
record the homology class of the curve in $X\setminus V$ rather than
in $X$.  These are additional data: two elements of $H_2(X\setminus V)$
represent the same element of $H_2(X)$ if they differ by an element of
the set ${\cal R}\subset H_2(X\setminus V)$ of rim tori (the name
refers to the fact that each such class can be represented by a torus
embedded in the boundary of a tubular neighborhood of $V$).  The
subtlety is that this homology information is intertwined with the
intersection data, and so the appropriate homology-intersection data form 
a covering space ${\cal H}^V_X$ of $H_2(X) \times V^\ell$ with fiber
${\cal R}$.  This is constructed in Section 5.

We then come to the key step of showing that the space $\M^V$ of
$V$-regular maps carries a fundamental homology class.  For this we
construct an orbifold compactification of $\M^V$ --- the space of
$V$-stable maps.  Since $\M^V$ is a union of open components of
different dimensions the appropriate compactification is obtained by
taking the closure of $\M_{g,n,s}^V(X,A)$ separately for each $g,n,A$
and $s$.  This is exactly the procedure one uses to decompose a
reducible variety into its irreducible components.  However, since we
are not in the algebraic category, this closure must be defined via
analysis.

The required analysis is carried out in Sections 6 and 7.  There we
study the sequences $(f_n)$ of $V$-regular maps using an iterated
renormalization procedure.  We show that each such sequence limits to a
stable map $f$ with additional structure. The basic point is that some
of the components of such limit maps have images lying in $V$, but
along each component in $V$ there is a section $\xi$ of the normal
bundle of $V$ satisfying an elliptic equation $D^N\xi=0$; this $\xi$
`remembers' the direction from which the image of that component came
as it approached $V$.  The components which carry these sections are
partially ordered according to the rate at which they approach $V$ as
$f_n\to f$.  We call the stable maps with this additional structure
`$V$-stable maps'.  For each $g,n,A$ and $s$ the $V$-stable maps form
a space $\ov{\M}_{g,n,s}^V(X,A)$ which compactifies the space of
$V$-regular maps by adding frontier strata of (real) codimension at
least two.

This last point requires that $(J,\nu)$ be $V$-compatible.  In Section
3 we show that for $V$-compatible $(J,\nu)$ the operator $D^N$
commutes with $J$.  Thus $\mbox{ker }D^N$, when nonzero, has (real)
dimension at least two.  This ultimately leads to the proof in Section
7 that the frontier of the space of $V$-stable maps has codimension at
least two.  In contrast, for generic $(J,\nu)$ the space of $V$-stable
maps is an orbifold with boundary and hence does not carry a
fundamental homology class.

The endgame is then straightforward.  The space of $V$-stable maps
comes with a map
\bear
\overline\M^V_{g,n,s}(X,A) \to \ov\M_{g,n+\ell(s)} \ti X^n \ti {\cal H}^V_X
\label{intro.2}
\eear
and relative invariants are  defined in exactly the same way
that the GW invariants are
defined from (\ref{intro.1}).  The new feature is the last factor, which
allows us to control how the images of the  maps intersect $V$.  Thus the
relative invariants give  counts of $V$-stable maps with constraints on the
complex structure of the domain, the
images of the marked points, and the  geometry of the intersection with $V$.

Section 1 describes the space of stable pseudo-holomorphic maps into a
symplectic manifold, including some needed features that are not yet
in the literature.  These are used in Section~2 to define the GW
invariants for symplectic manifolds and the associated invariants,
which we call GT invariants, that count possible disconnected curves.
We then bring in the symplectic submanifold $V$ and \pagebreak develop the ideas
described above.  Sections 3 and 4 begin with the definition of
$V$-compatible pairs and proceed to a description of the structure of
the space of $V$-regular maps.  Section 5 introduces rim tori and the
homology-intersection space ${\cal H}^V_X$.

For clarity, the construction of the space of $V$-stable maps is
separated into two parts.  Section 6 contains the analysis required
for several special cases with increasingly complicated limit maps.
The proofs of these cases establish all the analytic facts needed for
the general case while avoiding the notational burden of delineating
all ways that sequences of maps can degenerate. The key argument
is that of Proposition \ref{TaubesRenomlemma2}, which is
essentially a parametrized version of the original renormalization
argument of \cite{PW}. With this analysis in
hand, we define general $V$-stable maps in Section~7, prove the needed
tranversality results and give the general dimension count showing
that the frontier has sufficiently large codimension.  In Section~8
the relative invariants are defined and shown to depend only on the
isotopy class of the symplectic pair $(X,V)$.  The final section
presents three specific examples relating the relative invariants to
some standard invariants of algebraic geometry and symplectic
topology.  Further applications are  given in \cite{IP4}.

The results of this paper were announced in \cite{IP3}.  Related
results are being developed by by Eliashberg and Hofer \cite{E} and Li
and Ruan \cite{LR}.  Eliashberg and Hofer consider symplectic
manifolds with contact boundary and assume that the Reeb vector field
has finitely many simple closed orbits.  When our case is viewed from
that perspective, the contact manifold is the unit circle bundle of
the normal bundle of $V$ and all of its circle fibers -- infinitely
many -- are closed orbits.  In their first version, Li and Ruan also
began with contact manifolds, but the approach in the most recent
version of \cite{LR} is similar to that of \cite{IP3}. The relative
invariants we define in this paper are more general then those of
  \cite{LR} and appear, at least {\it a priori}, to give different gluing formulas.
 \def\sni#1{\smallbreak\noindent{#1}. }
\vglue12pt
\centerline{\bf Contents}

\sni{1}  Stable pseudo-holomorphic maps

\sni{2} Symplectic invariants

\sni{3} $V$-compatible perturbations

\sni{4} Spaces of $V$-regular maps

\sni{5} Intersection data and rim tori

\sni{6} Limits of $V$-regular maps

\sni{7} The space of $V$-stable maps

\sni{8} Relative invariants

\sni{9} Examples

\smallbreak\noindent\phantom{9. } Appendix


\section{Stable pseudo-holomorphic maps}

The  moduli space of $(J,\nu)$-holomorphic
maps from  genus $g$ curves  with\break $n$ marked points representing a class
$A\in H_2(X)$ has a compactification\break $\ov\M_{g,n}(X,A)$.  This
comes with a map
\begin{equation}
\ov\M_{g,n}(X,A)\ \longrightarrow \ \ov\M_{g,n}\ti X^n 
\label{modulimap1}
\end{equation}
where the first factor is the ``stabilization'' map st to the
Deligne-Mumford moduli space (defined by collapsing all unstable components
of the domain curve) and  the second factor records the images of the marked
points.  The compactification carries a `virtual fundamental class',
which, together with the map 
(\ref{modulimap1}), defines the Gromov-Witten invariants.

This picture is by now standard when $X$ is a K\"{a}hler manifold.
But in the general symplectic
case, the construction of  the compactification is scattered widely across
the literature (\cite{G}, \cite{PW},
\cite{P}, \cite{RT1},  \cite{RT2}, \cite{LT}, \cite{H} and \cite{IS})
and some needed
properties do not appear
explicitly anywhere.  Thus we
devote this section to reviewing and augmenting the construction of the
space of stable pseudo-holomorphic maps.

      Families of algebraic curves are well-understood from the work of Mumford
and others.  A smooth genus $g$ 
connected curve  $C$ with $n$ marked points is {\it stable} if
$2g+n\ge 3$, that is, if $C$ is
either a sphere with at least three marked
points, a torus with at least one marked point, or has genus $g\geq 2$. The
set of such curves, modulo
diffeomorphisms, forms the Deligne-Mumford moduli space $\M_{g,n}$.  This has
a compactification $\ov\M_{g,n}$ that is
a projective variety.  Elements of $\ov\M_{g,n}$ are called
`stable $(g,n)$-curves'; these are
unions of smooth stable components $C_i$ joined at $d$ double points with a
total of $n$ marked points and
Euler class $\chi(C)=2-2g+d$.  There is a universal curve
\bear
\ov\U_{g,n}=\ov\M_{g,n+1} \longrightarrow \ov\M_{g,n}
\label{1.univcurve}
\eear
whose fiber over each point of $[j]\in\ov\M_{g,n}$ is a stable curve $C$
in the equivalence class $[j]$ whenever $[j]$ has no automorphisms,
and in general is a curve $C/\Aut(C)$.  To avoid these quotients we
can lift to the moduli space of Prym structures as defined in
\cite{Lo}; this is a finite cover of the Deligne-Mumford
compactification and is a manifold. The corresponding universal curve
is a projective variety and is now a universal family, which we denote using the same
notation (\ref{1.univcurve}).  We also extend
this construction to the unstable range by taking  
$\ov\M_{0,n}=\ov\M_{0,3}$ for $n\le 2$ and 
$\ov\M_{1,0}=\ov\M_{1,1}$. We fix,
once and for all, a holomorphic embedding of $\ov\U_{g,n}$ into some $\P^N$.

At this juncture one has a choice of either working throughout with
curves with Prym structures, or working with ordinary curves and
resolving the orbifold singularities in the Deligne-Mumford space
whenever necessary by imposing Prym structures.  Moving between the
two viewpoints is straightforward; see Section~2 of \cite{RT2}.  To
keep the notation and discussion clear, we will consistently use
ordinary curves, leaving it to the reader to introduce Prym structures
when needed.

When one deals with maps $C\to X$ from a curve to another space one
should use a different notion of stability.  The next several
definitions define `stable holomorphic maps' and describe how they
form a moduli space.  We will use the term `special point' to refer to
a point that is either a marked point or a double point.

\numbereddemo{Definition} \label{defnstabledomain}
A {\it  bubble domain} $B$ of type $(g,n)$  is a finite connected union of
smooth oriented 2-manifolds
      $B_i$ joined at double points together with $n$ marked points, none of
which are double points. The  $B_i$, with their special points, are of
      two types:
\smallbreak
 (a)  stable components, and
\smallbreak  \hangindent=44pt\hangafter=1
 (b)  unstable rational components, called `\/{\it unstable bubbles}\/', which are
spheres with a complex structure and one or two special points.

\smallbreak\noindent
There must be at least one stable component.  Collapsing
the unstable
components to points gives  a {\it connected domain} ${\rm st}(B)$ which is a
stable genus~$g$ curve with $n$ marked points.
\enddemo

Bubble domains can be constructed from a stable curve by replacing
points by finite chains of 2-spheres.  Alternatively, they can be
obtained by pinching a set of
nonintersecting embedded circles (possibly contractible) in a smooth
2-manifold.  For our purposes, it
is  the latter viewpoint that is important.  It can be formalized as follows.

\numbereddemo{Definition}
A {\it resolution  of a $(g,n)$  bubble domain $B$ with $d$ double points}
is a  smooth oriented
2-manifold  with genus $g$, $d$ disjoint embedded circles $\gamma_{\ell}$,
and $n$ marked points
disjoint from the $\gamma_{\ell}$, together with a map `resolution map'
$$
r:\Si \to B
$$
that respects orientation and marked points, takes each
$\gamma_{\ell}$ to a double point of $B$, and restricts to a diffeomorphism
from the complement of the
$\gamma_{\ell}$ in $B$ to the complement of the double points.
\label{defnresolution}
\enddemo

We can put a complex structure $j$ on a bubble domain
$B$ by specifying  an
orientation-preserving  map
\begin{equation}
\phi_0: {\rm st}(B) \to \ov{\cal U}_{g,n}
\label{defnBdomainmap1}
\end{equation}
which is a diffeomorphism onto a fiber of $\ov\U_{g,n}$ and taking
$j=j_\phi$ to be $\phi^*j_{\cal U}$ on the stable components of $B$
and the standard complex structure on the unstable components.  We
will usually denote the complex curve $(B,j)$ by the letter $C$.

We next define $(J,\nu)$-holomorphic maps from bubble domains. These
depend on the choice of an $\w$-compatible almost complex structure
$J$ (see (\ref{A.1}) in the appendix), and on a `perturbation' $\nu$.
This $\nu$ is chosen from the space of sections of the bundle
$\mbox{Hom}(\pi^*_2T\P^N,\,\pi^*_1TX)$ over $X\ti \P^N$ that are
anti-J-linear:
$$
\nu(j_{P}(v))\ =\ -J(\nu(v))\qquad \forall v\in T\P^N
$$
where $j_{P}$ is the complex structure on $\P^N$. Let ${\cal J}$ denote the
space of such pairs $(J,\nu)$, and
fix one such pair.

\numbereddemo{Definition}
A $(J,\nu)$-\/{\it holomorphic map from a  bubble domain} $B$ is a  map
\begin{equation}
(f,\phi):B\ \longrightarrow\  X\times\ov{\cal U}_{g,n}\ \subset \  X\times\P^N
\label{2.defjholo}
\end{equation}
with  $\phi=\phi_0\circ {\rm st}$  as in (\ref{defnBdomainmap1})  such that, on
each component
$B_i$ of
$B$, $(f,\phi)$ is a smooth solution of the inhomogeneous Cauchy-Riemann
equation
\begin{equation}
\bar{\partial}_Jf\ =\ (f,\phi)^*\nu
\label{ CauchyRiemannEq}
\end{equation}
where $\bar{\partial}_J$ denotes the nonlinear elliptic operator
$\frac12(d+J_f\circ d\circ j_{\phi})$. In particular,
$\bar{\partial}_Jf=0$ on each unstable component.
\label{defnJnuholomorphic}
\enddemo

Each  map of the form (\ref{2.defjholo}) has   degree $(A,d)$ where
$A=[f(B)]\in H_2(X;\Z)$
and $d$ is the degree of
$\phi:{\rm st}(B)\to\P^N$;  $d\geq 0$ since $\phi$ preserves orientation and the
fibers of $\ov\U$ are holomorphic.   The ``symplectic area'' of the image is
the number
\begin{equation}
A(f,\phi)\ =\ \int_{(f,\phi)(B)} \w \times \w_{\P} \ =\ \int_{B} f^*\w
+\phi^* \w_{\P}\ =\ \w[A]+d \hskip.5in
\label{defA(f)}
\end{equation}
which depends only on the homology class of the map $(f,\phi)$.  Similarly,
the energy of  $(f,\phi)$ is
\begin{equation}
E(f,\phi)\ =\ \frac12\int_{B}|d\phi|_{\mu}^2+|df|_{J,\mu}^2\ d\mu\ =\
d +\frac12 \int_{B}|df|_{J,\mu}^2\ d\mu
\label{defEpfp}
\end{equation}
where $|\cdot|_{J,\mu}$ is the norm defined by the metric on $X$ determined
by $J$ and the metric $\mu$
on $\phi(B)\subset\P^N$.  These integrands are
conformally invariant, so
the energy depends only on $[j_{\phi}]$.  For  $(J,0)$-holomorphic maps,
the energy  and
the symplectic area are equal.

The following is the key definition for the entire theory.

\demo{Definition {\rm 1.4}}
A  $(J,\nu)$-holomorphic map $(f,\phi)$ is {\it stable} if each of its
component maps
$(f_i,\phi_i)=(f,\phi)|_{B_i}$ has positive energy.
\enddemo
\advance\theoremcount by 1

This means that each component $C_i$ of the domain is either a stable
curve, or else the image of $C_i$ carries a nontrivial homology class.

\proclaim{Lemma}\label{minenergylemma}
{\rm (a)}   Every $(J,\nu)$\/{\rm -}\/holomorphic  map has $E(f,\phi)\geq 1$.
\smallbreak

  {\rm (b)} There is a constant $0<\alpha_0<1${\rm ,} depending only on
$(X,J)${\rm ,} such that every component
$(f_i,\phi_i)$ of every stable $(J,\nu)$\/{\rm -}\/holomorphic  map into $X$
has $E(f_i,\phi_i) >\alpha_0$.
\smallbreak

 {\rm (c)} Every
$(J,\nu)$\/{\rm -}\/holomorphic  map $(f,\phi)$ representing a homology
class $A$  satisfies
$$
E(f,\phi)\leq \om(A)+C(3g-3+n)
$$
where  $C\ge 0$ is a constant which depends only on $\nu$ and the
metric on $X\ti \ov\U_{g,n}$ and which vanishes when $3g-3+n<0$. 
\endproclaim

\pf (a)\  If the  component maps  $(f_i,\phi_i)$ have   degrees
$(d_i,A_i)$ then\break $E(f,\phi)=\sum E(f_i,\phi_i)\geq \sum d_i$ by
(\ref{defEpfp}). But $\sum d_i\geq 1$ because at least one component is stable.

\smallbreak
(b)\  Siu and Yau \cite{SY} showed that there is a
constant $\alpha_0$, depending
only on $J$, such that
any smooth map $f:S^2\to X$ that is nontrivial in homotopy satisfies
$$
\frac12 \int_{S^2}|df|^2\ >\ \alpha_0.
$$
We may assume that $\alpha_0<1$.  Then  stable components have
$E(f_i,\phi_i)\geq 1$ as above, and each unstable
component either  has  $E(f_i,\phi_i)> \alpha_0$ or represents the trivial
homology class.  But in the latter case $f_i$ is  $(J,0)$-holomorphic, so
$E(f_i,\phi_i)=A(f_i,\phi_i)=\w[f_i]=0$, contrary to the definition of
stable map.
\smallbreak
(c) \ This follows from straightforward estimates using (\ref{
CauchyRiemannEq}) and (\ref{defEpfp}), and the observation that curves
in $\ov\M_{g,n}$ have at most $3g-3+n$ irreducible components. \enddemo

Let $\H^{J,\nu}_{g,n}(X,A)$ denote the set of $(J,\nu)$-holomorphic maps
from a  {\it smooth} oriented stable
Riemann surface with genus $g$ and $n$ marked points  to $X$ with
$[f]=A$ in $H_2(X;\Z)$.  Note that $\H$ is invariant under the group
$\mbox{Diff}(B)$ of
diffeomorphisms of the domain that preserve orientation and marked points:
if $(f,\phi)$ is
$(J,\nu)$-holomorphic then so is $(f\circ\psi,\phi\circ \psi)$ for any
diffeomorphism $\psi$.
Similarly, let $\ov\H^{J,\nu}_{g,n}(X,A)$ be the (larger) set of stable
$(J,\nu)$-holomorphic maps from a stable
$(g,n)$ bubble domain.

The main fact about $(J,\nu)$-holomorphic maps --- and the reason for
introducing bubble domains --- is the
following convergence theorem.  Roughly, it asserts that every sequence of
$(J,\nu)$-holomorphic maps from
a smooth domain has a  subsequence that converges modulo diffeomorphisms to a
stable map. This result,  first suggested by Gromov \cite{G},  is sometimes
called the
``Gromov Convergence Theorem''. The proof is the result of a series of
papers dealing with
progressively more general cases (\cite{PW}, \cite{P},
\cite{RT1}, \cite{H}, \cite{IS}).

\proclaimtitle{Bubble Convergence}
\proclaim{Theorem} \label{BTCthm}
Given any sequence $(f_j,\phi_j)$ of $(J_i,\nu_i)$\/{\rm -}\/holomorphic maps
with $n$ marked points{\rm ,} with
$E(f_j,\phi_j)<E_0$ and
$(J_i,\nu_i)\to(J,\nu)$ in $C^k${\rm ,} $k\ge 0${\rm ,} one can pass to a
subsequence  and find

\begin{itemize}
\item[{\rm (i)}] a  $(g,n)$ bubble domain $B$ with resolution  $r:\Si\to B${\rm ,} and

\item[{\rm (ii)}] diffeomorphisms $\psi_j$ of $\Sigma$ preserving the orientation and
the marked points{\rm ,}
\end{itemize}
so that the modified subsequence
$(f_j\circ\psi_j,\phi_j\circ \psi_j)$ converges to a limit
$$
\Si\stackrel{r}{\ \longrightarrow\ } B\stackrel{(f,\phi)}{\
\longrightarrow\ } X
$$
where  $(f,\phi)$ is a stable $(J,\nu)$\/{\rm -}\/holomorphic map.  This convergence
is in $C^0${\rm ,}  in $C^{k}$ on
compact sets not intersecting the collapsing curves $\gamma_{\ell}$ of the
resolution $r${\rm ,} and the area and
energy integrals {\rm (1.6)} and {\rm (\ref{defEpfp})} are preserved in the
limit.
\endproclaim

Under the convergence of Theorem \ref{BTCthm}, the image curves
$(f_j,\phi_j)(B_j)$ in $X\times \P^N$
converge to $(f,\phi)(B)$ in the Hausdorff distance $d_H$, and the marked
points and their images converge.
Define a pseudo-distance on  $\ov\H^{J,\nu}_{g,n}(X,A)$ by
\begin{eqnarray}
d\left( (f,\phi),(f',\phi')\right)& =& d_H\left(\phi(\Si),\phi'(\Si)\right) +
d_H\left(f(\Si),f'(\Si)\right)\label{Hausdorffdistance}
\\
&& +\ \sum
\ d_{X}\left(f(x_i),f'(x_i')\right) \nonumber
\end{eqnarray}
where the sum is over all the marked points $x_i$.  The {\it space of
stable maps}, denoted
$$
\ov{\M}^{J,\nu}_{g,n}(X,A) \qquad\mbox{or}\qquad \ov{\M}_{g,n}(X,A),
$$
is the space of equivalence classes in  $\ov\H^{J,\nu}_{g,n}(X,A)$, where two
elements are equivalent if the distance
(\ref{Hausdorffdistance}) between them is zero.  Thus  orbits of the
diffeomorphism group become single points in the quotient.  We always
assume the stability condition $2g+n\ge 3$.

The following structure theorem then follows from Theorem \ref{BTCthm}
above and the results of \cite{RT1} and
\cite{RT2}.  Its statement involves the canonical class $K_X$  of
$(X,\w)$  and the following two terms.

\numbereddemo{Definition} \label{def1.7}
(a)\   A symplectic manifold $(X,\w)$ is called {\it semipositive} if there is no spherical homology class $A\in
H_2(X)$ with $\omega(A)>0$ and  $0<2 K_X[A] \leq \dim X-6$.

\smallbreak

(b)\   A stable map $F=(f,\phi)$ is {\it irreducible} if it is
generically injective, i.e., if $F^{-1}(F(x))=x$ for generic points $x$.
\enddemo
 
    Let $\ov\M_{g,n}(X,A)^*$ be the moduli space of  irreducible stable
maps.  Definition (\ref{def1.7}b) is  equivalent to saying that  the
restriction of $f$ to the
union of the unstable components of its domain is generically injective
(such maps are called {\it simple} in {\rm \cite{MS}}).  Thus there
are two types
of reducible maps: \pagebreak maps whose restriction to some unstable rational
component factors through a covering map $S^2\to S^2$ of degree two or
more, and maps with two or more unstable rational components with
the same image.

\proclaimtitle{Stable Map Compactification} 
\proclaim{Theorem}  \label{BTStructurethm}
{\rm (a)}\  $\ov{\M}^{J,\nu}_{g,n}(X,A)$ is a compact metric space{\rm ,} and there
are continuous maps
\begin{equation}
\label{2.MainEvalmap}
\M^{J,\nu}_{g,n}(X,A)   \ \  \stackrel{\iota}{\hookrightarrow}  \ \
\ov{\M}^{J,\nu}_{g,n}(X,A) \ \ \lraop{{{\rm {st}} \times {\rm{ev}}}} \ \
\ov{\M}_{g,n}
\times X^n
\end{equation}
\noindent where $\iota$ is an embedding{\rm ,} ${\rm st}$ is the stabilization map
applied to the domain
$(B, j_\phi)${\rm ,} and
$\ev$  records the images of the marked points.  The composition
{\rm (\ref{2.MainEvalmap})} is smooth.

\smallskip

{\rm (b)}\    For generic $(J,\nu)$,
$\ov{\M}^{J,\nu}_{g,n}(X,A)^*$ is an oriented orbifold of {\rm (}\/real\/{\rm )} dimension
\begin{equation}
-2K_X[A]+(\dim X-6)(1-g)+2n.
\label{2.maindimcount}
\end{equation}
   Furthermore{\rm ,}  each stratum ${\cal S}_k^*\subset
\ov{\M}^{J,\nu}_{g,n}(X,A)^*$ consisting of maps whose domains have
$k$ double points is a suborbifold of
{\rm (}\/real\/{\rm )} codimension $2k$.

\smallskip

{\rm (c)}\    For generic $(J,\nu)${\rm ,} when $X$ is semipositive or
$\ov{\M}^{J,\nu}_{g,n}(X,A)$ is irreducible{\rm ,} then the image
     of $\ov{\M}^{J,\nu}_{g,n}(X,A)$ under ${\rm st}\ti {\rm ev}$
     carries a homology class.
   \endproclaim

The phrase `for generic $(J,\nu)$' means that the statement holds for all
$(J,\nu)$ in a second category subset  of the space (\ref{defJVeq}).

The manifold structure in (b) can be described as follows.  Given a
stable map $(f, \phi)$ with smooth domain $B$, choose a local
trivialization ${\cal U}_{g,n}=\M_{g,n}\times B$ of the universal
curve in a neighborhood $U$ of $\phi(B)$.  Then $\phi$ has the form
$([j_\phi],\psi)$ for some diffeomorphism $\psi$ of $B$, unique up to
$\mbox{Aut}(B)$ (and unique when $B$ has a Prym structure).  Then
\bear
{\cal S}_\phi\ =\ \left\{ \mbox{$(J,\nu)$-holomorphic}\
(f,\phi)\ |\ \phi=([j_\phi], \mbox{id.})\right\}
\label{1.slice}
\eear
is a slice for the action of the diffeomorphism group because any
$(f',\phi')= (f', [j_{\phi'}],\psi)$  with $\phi'(B)$ in $U$ is equivalent to
$(f'\circ\psi^{-1}, [j_{\phi'}], \mbox{id.})$, uniquely as above.  Thus the
space of stable maps is locally modeled by the product of
$\M_{g,n}$ and the set of $(J,\nu)$-holomorphic maps from the fibers
of the universal curve, which is a  manifold as in \cite{RT2}.

The strata ${\cal S}_k^*$ are orbifolds because with  irreducible maps one
can use variations in the pair $(J,\nu)$ to achieve the
tranversality needed to show that the moduli space
is locally smooth  and oriented for generic $(J,\nu)$.   This is
proved in  Lemma 4.9
in \cite{RT1} and Theorem 3.11 in \cite{RT2} (the proof also applies to irreducible maps with  ghost bubbles, which are
unnecessarily singled out in  \cite{RT1}). Moreover, the gluing
theorem of Section 6 of \cite{RT1} proves
that  ${\cal S}_k^*$ has an  orbifold tubular neighborhood in
$\ov{\M}^{J,\nu}_{g,n}(X,A)^*$.

Theorem \ref{BTStructurethm}c was proved in \cite{RT2} for
semipositive $(X,\w)$ by reducing the moduli space as follows. Every
reducible stable map $f\in \ov\M_{g,n}(X,A)$ factors through an
irreducible stable map $f_0\in \ov\M_{g,n}(X,A_0)^*$ which has the
same image as $f$, with the homology classes satisfying
$\om([f_0])\le \om ([f])$.  Replacing each reducible $f$ by $f_0$
yields a `reduced moduli space' without reducible maps   whose image
under ${\rm st}\ti {\rm ev}$ contains the image of the original moduli space.
Semipositivity then implies that all boundary strata of the image of
the reduced moduli space are of codimension at least {\rm 2}.

\demo{{R}emark {\rm 1.9 (Stabilization)}} The semipositive assumption
in Theorem \ref{BTStructurethm}c can be removed
in several ways (\cite{LT}, \cite{S}, \cite{FO}, \cite{R}), each  leading
to  a moduli space which carries a  ``virtual fundamental class'', or
at least whose image defines a homology
class as in Theorem \ref{BTStructurethm}c. Unfortunately these
approaches involve replacing the space of $(J,\nu)$-holomorphic maps
with a more complicated and  abstract space.  It is preferable, when
possible, to work directly with  $(J,\nu)$-holomorphic maps
where one can use the equation (\ref{ CauchyRiemannEq}) to make
specific geometric and P.D.E.\ arguments.

In a separate paper \cite{IP5} we describe an alternative approach
based on the idea of adding enough additional structure to insure that
all stable $(J,\nu)$-holomorphic maps are irreducible.  More
specifically, we develop a scheme for constructing a new moduli space
$\tilde{\M}$ by consistently adding additional marked points to the
domains and imposing constraints on them in such a way that (i) all
maps in $\tilde{\M}$ are irreducible, and (ii) $\tilde{\M}$ is a
finite (ramified) cover of the original moduli space.  Theorem
\ref{BTStructurethm}c then applies to $\tilde{\M}$ and hence
$\tilde{\M}$, divided by the degree of the cover, defines a homology
class.  \label{1.remark}
\enddemo

\section{Symplectic invariants}

For generic $(J,\nu)$ the space of
stable maps carries a
fundamental homology class.   For each $g,n$ and $A$, the pushforward of
that class under the evaluation
map (\ref{modulimap1}) or (\ref{2.MainEvalmap}) is the `Gromov-Witten'
homology class
\begin{equation}
\left[\ov \M_{g,n}(X,A)\right] \ \in\
H_*(\ov\M_{g,n};\Q)\otimes H_*(X^n;\Q).
\label{modulimap2}
\end{equation}
A cobordism argument shows that this is independent of the choice of generic
$(J,\nu)$, and hence depends
only on the symplectic manifold $(X,\w)$.  Frequently, this Gromov-Witten
invariant is thought of  as the
collection of numbers obtained by evaluating (\ref{modulimap2}) on a basis
of the dual cohomology group.

For our purposes it is convenient to assemble the GW invariants into power
series in such a way that
disjoint unions of maps correspond to products of the power series.  We define those series in this section.  Along
the way we describe the geometric interpretation of the invariants.

Let $NH_2(X)$ denote the Novikov ring as in \cite{MS}.  The elements
of $NH_2(X)$ are sums $\sum c_A t_A$ over $A\in H_2(X;\Z)$ where
$c_A\in\Q$, the $t_A$ are variables satisfying $t_At_{B}=t_{A+B}$,
$c_A=0$ if $\omega(A)<0$, and where, for each $C>0$ there are only 
finitely many nonzero coefficients $c_A$ with energy $\w(A)\leq C$.  
After summing on $A$ and dualizing, (\ref{modulimap2}) defines a map
\begin{equation}
{\rm GW}_{g,n}:\, H^*(\overline{\M}_{g,n})\otimes H^*(X^n) \ \to \ NH_2(X).
\label{1.3}
\end{equation}

We can also sum over $n$ and $g$ by setting $\overline{\M} =
\bigcup_{g,n}\overline{\M}_{g,n}$, letting ${\Bbb T}^*(X)$ denote the
total (super)-tensor algebra ${\Bbb T}(H^*(X))$ on the rational
cohomology of $X$, and introducing a variable $\lambda$ to keep track
of the Euler class.  The total Gromov-Witten invariant of $(X,\w)$ is
then the map
\bear
  {\rm GW}_X: H^*(\ov\M)\otimes {\Bbb T}^*(X)\ra NH_2(X)[\la].
\label{defGW}
\eear
defined by the Laurent series
\begin{equation}
{\rm GW}_{X}\ =\ \sum_{A,g,n} \frac{1}{n!}\,{\rm GW}_{X,A,g,n}\ t_A\ \la^{2g-2}.
\label{defn.7}
\end{equation}
The diagonal action of the symmetric group $S_n$ on
$\overline{\M}_{g,n}\ti X^n$ leaves ${\rm GW}_X$ invariant up to sign, and if
$\kappa\in H^*(\ov\M_{g,n})$ then
${\rm GW}_X(\kappa, \alpha)$ vanishes unless $\alpha$ is a tensor of length $n$.

We can recover the familiar geometric interpretation of these invariants
by  evaluating on cohomology classes.  Given $\kappa\in
H^*(\ov\M;\Q)$ and  a
vector $\alpha=(\alpha_1,\dots, \alpha_n)$ of rational cohomology classes
in $X$ of length $n=\ell(\alpha)$,
fix a generic $(J,\nu)$ and generic
geometric representatives $K$ and $A_i$ of the Poincar\'e duals of
$\kappa$ and of the $\alpha_i$ respectively. Then
${\rm GW}_{X,A,g,n}(\kappa,\alpha)$ counts, with orientation,  the number of genus
$g$ $(J,\nu)$-holomorphic maps
$f:C\to X$ with $C\in K$ and  $f(x_i)\in A_i$ for each of the $n$ marked
points $x_i$.
By the usual dimension counts, this
vanishes   unless
$$
\deg \kappa+\sum\deg \alpha_i-2\ell(\alpha)=
(\mbox{dim}\, X-6)(1-g)-2K_X[A].
$$

It is sometimes useful  to incorporate the so-called  `$\psi$-classes'.   There
are canonically oriented real 2-plane bundles ${\cal L}_i$ over
$\overline{\M}_{g,n}(X,A)$
whose fiber at each map $f$ is the  cotangent space to the (unstabilized)
domain curve
at the $i^{\rm th}$ marked point.  Let $\psi_i$ be the Euler class of
${\cal L}_i$, and for each vector $D=(d_1,\dots d_n)$ of nonnegative
integers let $\psi_D=\psi_1^{d_1}\cup\dots \cup\psi_n^{d_n}$.  Replacing
the left-hand side of (\ref{modulimap2}) by the pushforward of the cap product
$\psi_D\cap \left[\ov\M_{g,n}(X,A)\right]$
and again dualizing gives invariants
\begin{equation}
{\rm GW}_{X,g,n,D}:\, H^*(\ov\M_{g,n})\otimes H^*(X^n) \ \to \ NH_2(X) \pagebreak
\label{1.45}
\end{equation}

\noindent
which agree with (\ref{1.3}) when $D$ is the zero vector.  These invariants
can be included in ${\rm GW}_X$ by adding variables in the series (\ref{defn.7})
which keep
track of the vector $D$. To keep the notation manageable we will leave that
embellishment to the reader.

The  GW invariant (\ref{defGW}) counts
$(J, \nu)$-holomorphic maps from {\it connected}  domains.  It is often
more natural to work with maps whose domains are disjoint unions.  Such
$J$-holomorphic curves arose, for
example, in Taubes' work on  the Seiberg-Witten invariants (\cite{T}).
In fact, there is a  simple and natural way of
extending  (\ref{defGW}) to this more general case.

       Let $\wt{\M}_{\chi,n}$
be the space of all compact Riemann surfaces of Euler characteristic $\chi$
with finitely many unordered components  and with a total of
$n$ (ordered) marked points.  For each such surface, after we fix an ordering
of its components, the locations of the marked
points define  an ordered partition
$\pi=(\pi_1,\dots,\pi_l)\in{\cal P}_n$. Hence
\best
\wt \M_{\chi,n}=\ma\bigsqcup_{\pi\in {\cal P}_n}\;
\ma\bigsqcup_{g_i}\  \l(\;\ov\M_{g_1,\pi_1}\ti \dots \ti
      \ov\M_{g_l,\pi_l}\r)/S_l
\eest
where ${\cal P}_n$ is the set of all ordered partitions of the set
$\{x_1,\dots x_n\}$, $\ov\M_{g_i,\pi_i}$ is the  space of stable curves
with $n_i$ marked points labeled by   $\pi_i$, and where the second union
is over all $g_i$ with $\sum(2-2g_i)=\chi$. The symmetric group $S_l$
acts by interchanging the components.  Define the  ``Gromov-Taubes''
invariant
\bear
{\rm GT}_X: H^*(\wt{\M})\otimes {\Bbb T}^*(X) \ \to \ NH_2(X)[\la]
\label{defnGT-1}
\eear
by
\begin{equation}
{\rm GT}_X\ =\ e^{{\rm GW}_X}.
\label{defnGT}
\end{equation}
This exponential uses the ring structure on both sides of
(\ref{defnGT-1}).  Thus  for $\al=\al_1\otimes\dots \otimes \al_n$ and
      $\kappa=\kappa_1\otimes \dots \otimes \kappa_l$,  
\best
&& \hskip-24pt
{\rm GT}_{X, n}(\kappa,\alpha)\\ &&\quad =\;\sum_{\pi\in{\cal P}_n}
{\ep(\pi) \over l! }{ n \choose n_1, \dots, n_l}
      {\rm GW}_{X,n_1}(\kappa_{1},\al_{\pi_1})
\otimes  \dots \otimes {\rm GW}_{X, n_l}(\kappa_{l},\alpha_{\pi_l})
\eest
where, for each partition $\pi=(\pi_1,\dots \pi_l)$,  $\al_{\pi_i}$ is
the product of $\al_j$ for all $j\in \pi_i$
and $\ep(\pi)=\pm 1$ depending on the sign of the permutation
$(\pi_1,\dots, \pi_l)$ and the degrees of $\al$.

As before, when (\ref{defnGT}) is expanded as a Laurent series,
\vglue12pt
\centerline{${\displaystyle
{\rm GT}_X(\kappa, \alpha)\ =\ \sum_{A,\chi,n}
      {1\over n!}\;{\rm GT}_{X,A,\chi,n}(\kappa,\alpha)\;t_A \ \la^{-\chi},
}$} \pagebreak\noindent
the coefficients count the number of curves (not necessarily connected)
with Euler characteristic $\chi$ representing $A$ satisfying the
constraints $(\kappa,\alpha)$.  Note that  $A$ and $\chi$ add when one
takes disjoint unions, so that the variables $t_A$ and $\la$ multiply.

 \section{$V$-compatible perturbations}

We now begin our main task: extending the symplectic invariants of
Section 2 to invariants of $(X,\w)$ {\it relative to a codimension two
symplectic submanifold $V$}. Curves in $X$ in general position will
intersect such a submanifold $V$ in a finite collection of points.
Our relative invariants will still be a count of $(J,\nu)$-holomorphic
curves in $X$, but will also keep track of how those curves intersect
$V$.  But, instead of generic $(J,\nu)$, they will count holomorphic
curves for special $(J,\nu)$: those `compatible' to $V$ in the sense
of Definition \ref{def3.2} below.

Because $(J,\nu)$ is no longer generic, the construction of the space
of stable maps must be thought through again and modified.  That will
be done over the next six sections.  We begin in this section by
developing some of the analytic tools that will be needed later.

       The universal moduli space of stable maps
$\ov{\UM}_{g,n}(X)\to \J$ is the set
of all maps into $X$ from {\it some} stable $(g,n)$ curves which are
$(J,\nu)$-holomorphic for {\it some} $(J,\nu)\in {\cal J}$.  If we fix a
genus $g$ two-manifold $\Si$, this is the set of
$(f,\phi, J,\nu)$ in $\mbox{Maps}(\Si, X\ti \ov{\cal U}_{g,n}) \ti \J$ with
$\bar\partial_{J}f=\nu$. Equivalently, $\ov{\UM}_{g,n}(X)$ is
the  zero set of
\bear
\Phi(f,\phi,J,\nu)\ = \frac12\left(\ df+J\circ df\circ j\right)\,-\,\nu
\label{defofPhi}
\eear
where $j$ is the complex structure  on the domain determined by  $\phi$.
We will often abuse notation by writing $j$ instead of $\phi$.

In a neighborhood of $(f,\phi)$ the space of stable maps  is modeled by
the slice (\ref{1.slice}).  Within that slice, the variation in
$\phi$ lies in the tangent space to $\ov\M_{g,n}$, which is canonically
identified with $H^{0,1}(TC)$ where $C$ is the image of~$\phi$.

\proclaim{Lemma}
\label{linearizationlemma}
The linearization of {\rm (\ref{defofPhi})}  at a point $(f,j,J,\nu)\in \UM_{g,n}$
is the
elliptic operator
$$
D\Phi: \Gamma(f^*TX)\oplus H^{0,1}_j (TC)
\oplus \mbox{\rm End}(TX,J)\oplus {\rm Hom}_J(T\P^N, TX)\to
\Omega^{0,1}_j(f^*TX)
$$
given by
$$
D\Phi(\xi,k,K,\mu) =D_f (\xi,k) + {1\over 2}K f_* j - \mu %
$$
where $C$ is the domain of $f$ and   
$D_f(\xi,k)=D\Phi(\xi,k,0,0)$ is \pagebreak defined by

\begin{eqnarray}
\label{linearizationequation}
&&\nonumber\\
\noalign{\vskip-32pt}
D_f(\xi,k)(w)& =& \frac12\left[ \nabla_w \xi+J\nabla_{jw}\xi +
(\nabla_\xi J)(f_*(jw))+ J f_* k(w) \right]\\
&&-\ (\nabla_\xi \nu)(w) \nonumber
\end{eqnarray}
for each  vector $w$ tangent to the domain{\rm ,} where $\nabla$ is the
pullback connection on
$f^*TX$.
\endproclaim

\pf The variations with respect to $j,J$ and $\nu$ are obvious (cf.\
equation (3.9) in \cite{RT2}), so we need
only check the variation with respect to $f$.
The calculation in \cite[Lemma 6.3]{RT1} gives
\begin{eqnarray*}
D_f (\xi,k)(w)& =&
      \frac12\left[ \nabla_w\xi+J\nabla_{jw}\xi+ \frac12
(\nabla_\xi J)(f_*(jw)+Jf_*(w)) +J f_* k(w)  \right] \\
&&-\ (\nabla^J_\xi \nu)(w)
\end{eqnarray*}
where $\nabla^J=\nabla +\frac12 (\nabla J)J$. By the equation
$\Phi(f,j,J,\nu)=0$, this
agrees with (\ref{linearizationequation}).
\enddemo

As mentioned above, we will restrict attention to a
subspace of ${\cal J}^V$ consisting of pairs  $(J,\nu)$ that are compatible
with $V$ in the following
sense.  Denote the orthogonal projection onto the normal bundle
$N_V$ by $\xi\mapsto\xi^N$; this uses the metric defined by $\omega$
and $J$ and hence depends on $J$.

\numbereddemo{Definition}
\label{def3.2}
Let ${\cal J}^V$ be the submanifold of ${\cal J}$ consisting of pairs
$(J,\nu)$ whose 1-jet along $V$ satisfies the following  three conditions:
$$
\hskip-1.2in {\rm (a)}\ \ \    \mbox{ $J$ preserves $TV$  and
$\nu^N|_V=0$,}\\
$$
\noindent and for all   $\xi\in N_V$, $v\in TV$ and $w\in TC$
      \begin{eqnarray}
\label{defJVeq}
      &\mbox{(b)} & \ \
\left[\l(\nabla_{\xi}J+J\nabla_{J\xi}J\r)(v)\right]^N\ =\
\left[(\nabla_v
J)\xi+J(\nabla_{Jv}J)\xi\right]^N, \\
      &\mbox{(c)} & \ \
\left[\l(\nabla_{\xi}\nu+J\nabla_{J\xi}\nu\r)(w)\right]^N\ =\
\left[(J\,\nabla_{\nu(w)}J)\xi \right]^N.\nonumber
\end{eqnarray}
\enddemo

The first condition means that $V$ is a $J$-holomorphic submanifold, and
that $(J,\nu)$-holomorphic
curves in $V$ are also $(J,\nu)$-holomorphic in $X$.   Conditions (b) and (c)
relate to the variation of
such maps; they are chosen to ensure that Lemma \ref{DNcxlemma} below holds.
Condition (b) is equivalent to the vanishing of some of the components of
the Nijenhuis tensor $N_J$ along $V$, namely that the normal component of
$N_J(v,\xi)$ vanishes whenever $v$ is tangent and $\xi$ is normal to $V$.  Thus
(b)  can be thought of as the `partial integrability' of $J$ along $V$.

For each  $(J,\nu)$-holomorphic map $f$ whose image lies in $V$, we obtain
an operator  $D_f^N: \Gamma(f^*N_V)\to \Omega^{0,1}(f^*N_V)$ by restricting
      the linearization (\ref{linearizationequation}) to the normal bundle:
\bear
D^N_f(\xi)\ =\ \left[D_{f}(\xi,0)\right]^N.
\label{3.4}
\eear
\proclaim{Lemma}
\label{DNcxlemma}
Let $(J,\nu)\in {\cal J}^V$.  Then for each $(J,\nu)$\/{\rm -}\/holomorphic
map $f$ whose image
lies in $V${\rm ,}   $D^N_f$ is a complex operator
{\rm (}\/that is{\rm ,} it commutes\break with $J${\rm ).}
\endproclaim

\pf   Since $J$ preserves the normal bundle, we must verify that 
$[D(J\xi)-JD(\xi)]^N=0$ for each
$\xi\in N_V$.  By (\ref{linearizationequation}), the quantity
$2J[D(J\xi)-JD(\xi)](w)$ is
$$
(J\nabla_w J)\xi -(\nabla_{jw}J)\xi+{1\over 2}[\nabla_\xi
J+J\nabla_{J\xi}J](f_*(jw)+Jf_*(w))-2[\nabla_\xi\nu
+J\nabla_{J\xi}\nu](w).
$$
After substituting $f_*(w)=2\nu(w)-Jf_*(jw)$ into the first term  and
writing $v=f_*(jw)$, this becomes
\begin{eqnarray*}
 2(J\nabla_{\nu(w)} J)\xi&\hskip-8pt -\hskip-8pt& (J\nabla_{Jv} J)\xi -(\nabla_{v}J)\xi
+\nabla_\xi J(v)\\
& \hskip-8pt+\hskip-8pt&J\nabla_{J\xi}J(v)-2\nabla_\xi\nu(w) -2J\nabla_{J\xi}\nu(w).
\end{eqnarray*}
Taking the normal component, we see that the sum of the second, third, fourth, and
fifth  terms
vanishes by (\ref{defJVeq}b), while the sum of the first, sixth, and
seventh terms vanishes by
(\ref{defJVeq}c). \enddemo

We conclude this section by giving a local normal form for holomorphic maps
near the points where they intersect $V$.  This will be used repeatedly
later.  The proof  is adapted from McDuff \cite{M}.

Here is the context.   Let $V$ be a codimension two $J$-holomorphic
submanifold of $X$ and $\nu$ be a perturbation that vanishes in the
normal direction to $V$
as in (\ref{defJVeq}).  Fix  a local holomorphic coordinate
$z$ on an open set ${\cal O}_{C}$ in a Riemann surface $C$.
Also fix local coordinates $\{v^i\}$ in an open set ${\cal O}_V$ in
$V$ and extend these to local coordinates $(v^i,x)$ for $X$ with
$x\equiv 0$ along $V$ and so that  $x=x^1+ix^2$ along $V$  with
$J(\partial/\partial x^1)=\partial/\partial x^2$ and
$J(\partial/\partial x^2)=-\partial/\partial x^1$.

\proclaimtitle{normal form}  \proclaim{Lemma} \label{ContactLemma}
 Suppose that $C$ is a smooth connected curve
and  $f: C \rightarrow X$ is a
$(J,\nu)$\/{\rm -}\/holomorphic map that intersects $V$  at a point $p=f(z_0) \in
V$ with $z_0\in {\cal O}_{C}$ and $p\in {\cal O}_V$.  Then either
{\rm (i)} $f(C) \subset V${\rm ,} or
{\rm (ii)}  there is an  integer $d>0$ and a nonzero $a_0\in\cx$  so
that in the above coordinates
\begin{equation}
f(z,\bar{z})\ =\ \left(\,p^i+O(|z|),\  a_0z^d+O(|z|^{d+1})\,\right)
\label{Taylorexpansion}
\end{equation}
where $O(|z|^k)$ denotes a function of $z$ and $\bar{z}$ that vanishes to order
$k$ at $z=0$.
\endproclaim

\pf   Let $J_0$ be the standard complex structure in
the coordinates $(v^i,x^{\alpha})$.   The components of the matrix of
$J$ then satisfy
\begin{equation}
(J-J_0)^i_j  = O(|v|+|x|), \qquad (J-J_0)^{\alpha}_{\b} = O(|x|),
\qquad
(J-J_0)^i_{\alpha}  = O(|x|).
\label{AronLemma1}
\end{equation}
Set
$$ A = (1 - J_0 J)^{-1}  (1+J_0 J) \qquad \mbox{and} \qquad
\hat{\nu} = 2(1 - J J_0)^{-1}\nu.
$$
With the usual definitions $\bar{\partial}  f = \frac{1}{2} (df +
J_0 dfj)$ and $\partial f = \frac{1}{2}(df - J_0 dfj)$, the $(J,\nu)$-holomorphic map equation $\bar{\partial}_Jf =
\nu$ is equivalent to
\begin{equation}
\bar{\partial} f = A \partial f + \hat{\nu}.
\label{AronLemma2}
\end{equation}
Conditions (\ref{AronLemma1})  and the fact that the normal component of
$\nu$ also vanishes along $V$ give
$$
A^i_j = O(|v|+|x|), \qquad A^{\alpha}_{\b} = O(|x|), \qquad
A^i_{\alpha} = O(|x|),
\qquad \nu^{\alpha} = O(|x|).
$$
Now write $f=(v^i(z,\bar{z}),x^\alpha(z,\bar{z}))$.    Because
$A^i_{\alpha}$ vanishes along $V$ and the functions $|dv^i|$  and $\partial
A^i_{\alpha} /\partial x^{\b}$
are bounded  near $z_0$, we obtain
$$
\left|dA^i_{\alpha}\right| \ \leq\  \left|\frac{\partial
A^i_{\alpha}}{\partial v^j}
\cdot d v^j \   + \
\frac{\partial A^i_{\alpha}}{\partial x^{\b}} \cdot dx^{\b} \right|
\ \leq \ c\left(|x| + |dx|\right).
$$
Since  $\nu^{\alpha}$ also vanishes along $V$ by Definition \ref{defJVeq}a,
we get exactly the same bound on
$|d\nu^{\alpha}|$.  Returning to  equation (\ref{AronLemma2}) and looking
at the $x$ components, we have
\begin{equation}
\bar{\partial} x^{\alpha} = A^{\alpha}_i\partial v^i + A^{\alpha}_{\b}
\partial x^{\b}+
\hat{\nu}^{\alpha},
\label{AronLemma3}
\end{equation}
and hence
$$
\partial \bar{\partial} x^{\alpha}\ =\  \partial A^{\alpha}_i \ \partial
v^i + A^{\alpha}_i \ \partial^2
v^i + \partial A^{\alpha}_{\b} \ \partial x^{\b} + A^{\alpha}_{\b}\
\partial^2 x^{\b}
+\partial\hat{\nu}^{\alpha}.
$$
Because $ \partial \bar{\partial} x^{\alpha}=2\Delta x^{\alpha}$ and the
derivatives of $v$ and $x$ are
locally bounded this gives
$$ |\Delta x^{\alpha}|^2\  \leq\  c \left(|x|^2 + |\partial x|^2\right).
$$
If $x^{\alpha}$ vanishes to infinite order at $z_0$ then Aronszajn's Unique
Continuation theorem (\cite[Remark 3]{A})  implies that
$x^{\alpha}\equiv 0$  in a neighborhood of $z_0$, i.e. $f(C)\subset V$
locally.  This
statement is independent of coordinates.  Consequently, the set of $z \in
C$ where $f(z)$ contacts $V$
to infinite order is both open and closed, so that
$f(C)\subset V.$  On the other hand, if the order of vanishing is finite, then
$x^{\alpha}(z,\bar{z})$ has a
Taylor expansion beginning with
$\sum_{k=0}^d a_k \bar{z}^kz^{d-k}$ for some $0<d<\infty$.  Since
$A^{\alpha}_{i}, A^{\alpha}_{\beta}$ and
$\nu^{\alpha}$ are all $O(|x|)$ and $x$ is $O(|z|^d)$, (\ref{AronLemma3}) gives
$$
\bar{\partial}x^{\alpha}\ =\ O(|x|)\ =\ O(|z|^d).
$$
Differentiating, we conclude that the leading term is simply $a_0z^d$.
This\break gives~(\ref{Taylorexpansion}).  \enddemo
 
\vglue-18pt

 \section{Spaces of $V$-regular maps}

\vglue-6pt

We have chosen to work with holomorphic maps for $(J,\nu)$ compatible with
$V$.  For these special
$(J,\nu)$ one can expect  more holomorphic curves than are present for a
completely general  choice of $(J,\nu)$.
In particular, with our choice, any $(J,\nu)$-holomorphic map
into $V$ is automatically holomorphic as a map into $X$.  Thus we have
allowed stable holomorphic maps that are
badly nontransverse to $V$ --- entire components can be mapped into $V$.
We will  exclude such maps and define
the  relative invariant  using only    `$V$-regular' maps.

\numbereddemo{Definition}
\label{defn4.Vregular}
A stable $(J,\nu)$-holomorphic map into $X$  is called\break {\it $V$\/{\rm -}\/regular}
if  no component of its domain is mapped entirely into $V$ and if none of
the special points (i.e. marked or double points) on its domain are
mapped into $V$.
\enddemo

The $V$-regular maps (including those with nodal domain) form an open subset
of the space of stable maps, which we denote by $\M^V(X,A)$. In this section
we will show how $\M^V(X,A)$ is  a disjoint union of components, and how the 
irreducible part of each component is an orbifold for generic $(J,\nu)\in{\cal J}^V$.

     Lemma \ref{ContactLemma} tells us that for each $V$-regular
map $f$,  the inverse image
$f^{-1}(V)$  consists of isolated  points $p_i$ on the domain $C$
distinct from the special
points.  It also shows that each $p_i$ has
a well-defined multiplicity $s_i$ equal to the order of contact of the
image of $f$ with $V$ at $p_i$.
The list of  multiplicities  is a vector $s=(s_1, s_2,\dots, s_\ell)$ of
integers $s_i\geq 1$.  Let ${\cal S}$ be the
set of all such vectors and define the degree, length, and order of
$s\in{\cal S}$ by
$$
\mbox{deg }s =\sum s_i ,\qquad \ell(s) = \ell, \qquad
|s|=s_1 s_2\cdots s_\ell.
$$
These vectors $s$ label the components of $\M^V(X,A)$:  associated
to each $s$ such that $\deg s=A\cdot V$ is
the  space
$$
\M_{g,n,s}^V(X,A)\,\subset\ \M_{g,n+\ell(s)}(X,A)
$$
of all $V$-regular maps $f$ such that $f^{-1}(V)$ is exactly the
marked points $p_i$, $1\leq i \leq \ell(s)$, each with multiplicity $s_i$.
Forgetting these last $\ell(s)$ points defines a projection
\bear
\begin{array}{c}
\M_{g,n,s}^V(X,A) \\[5pt]
\big\downarrow\\[5pt]
\M_{g,n}^V(X,A)
\end{array}
\label{ellfoldcover}
\eear
onto one  component of $\M_{g,n}^V(X,A)$, which is the
disjoint union of such
components.  Notice that for each $s$ (\ref{ellfoldcover}) is a
covering space whose group of deck
transformations is the group of renumberings of the last $\ell(s)$ marked
points.

\proclaim{Lemma} \label{Lemma5.1}
 For generic $(J,\nu)${\rm ,} the irreducible part of
$\M_{g,n,s}^V(X,A)$ is an orbifold with
\begin{eqnarray}\label{5.1}
{\rm dim}\  \M_{g,n,s}^V(X,A)& \hskip-6pt=\hskip-6pt& -2K_X[A]+({\rm dim }\,
X-6)(1-g)\\
&\hskip-6pt\hskip-6pt&+\ 2(n+\ell(s) -\deg s).\nonumber
\end{eqnarray}
\endproclaim
 
\pf We need only to show that the universal moduli space
$\UM_{g,n,s}^*$ is a manifold (after passing to Prym covers); the  Sard-Smale
theorem then implies that for generic $(J,\nu)$ the  moduli space
$\M_{g,n,s}^V(X,A)^*$ is an orbifold of dimension equal to
the (real) index of the linearization, which is precisely
(\ref{5.1}).

First, let ${\cal F}^V_{g,n}$ be the space of all data $(J,\nu,f,j,
x_1, \dots, x_n)$ as in \cite[Eq.\ (3.3)]{RT2}, but now taking $f$ to be
$V$-regular and $(J,\nu)\in {\cal J}^V$.  Define $\Phi$ on ${\cal
F}^V_{g,n}$ by $\Phi(J,\nu,f,j,\{x_i\})=\overline{\partial}_{jJ}f-\nu$.  The
linearization $D\Phi$ is onto exactly as in equations (3.10) and
(3.12) of \cite{RT2}, so that the universal moduli space
$\UM^{V*}_{g,n}=\Phi^{-1}(0)$ is smooth and its dimension is given by
(\ref{5.1}) without the final `$s$' terms.

It remains to show that the contact condition corresponding to each
ordered sequence $s$ is transverse; that will imply that $\UM_{g,n,s}^V(X,A)^*$ is a
manifold.  Consider the space $\div^d(C)$ of degree $d$ effective
divisors on $C$. This is a smooth manifold of complex dimension $d$.
(Its differentiable structure is as described in \cite[p.~236]{GH}:
given a divisor $D_0$, choose local holomorphic coordinates $z_k$
around the points of $D_0$; nearby divisors can be realized as the
zeros of monic polynomials in these $z_k$ and the coefficients of
these polynomials provide a local chart on $\mbox{Div} ^d(C)$.)
Moreover, for each sequence $s$ of degree $d$, let $\div_s(C)\subset
\div^d(C)$ be the subset consisting of divisors of the form $\ma\sum
s_ky_k$. This is a smooth manifold of complex dimension $\ell(s)$.

For each  sequence $s$ of degree $d$  define a map
$$
\Psi_s \ :\ \UM^V_{g,n+\ell(s)}\longrightarrow \div ^d(C)\ti
\div_s(C)
$$
by
$$ \Psi_s(J,\nu,f,j,\{x_i\},\{y_k\})\ =\  \left(f^{-1}(V),
\ \sum_k s_ky_k\right)
$$
where the $y_k$ are the last $\ell(s)$ marked points.
By Lemma \ref{ContactLemma}, there are local coordinates $z_l$ around the
points $p_l\in C$ and $f(p_l)\in V\subset X$ such that the leading term of
the normal component of $f$ is $z_l^{d_l}$;
hence
\best
\Psi_s(J,\nu,f,j,\{x_i\},\{y_k\})\ =\ \left(\ma\sum d_l p_l,
\ma\sum s_k y_k\right)
\eest
with $d_l\ge 1$, $\sum d_l=A\cdot V=d$. Let $\Delta\subset \div ^d(C)\ti
\div_s(C)$ denote the  
diagonal of $\div_s(C)\times \div_s(C)$. Then
\bear
\UM^V_{g,n,s}=\Psi_s^{-1}( \Delta).
\label{4.defUMs}
\eear
This is a manifold provided that $\Psi_s$ is transverse to $\Delta$.
Thus it suffices to show that at each fixed
$(J,\nu,f,j,\{x_i\},\{y_k\})\in \UM^V_{g,n,s}$ the differential $D\Psi$ is
onto the tangent space of the first factor.

To verify that, we need only to construct a deformation
$$(J,\nu_t,f_t,j,\{x_i\},\{y_k\})$$ that is tangent to
$\UM_{g,n+\ell(s)}$ to first order in $t$, where the zeros of
$f^N_t$ are, to first order in $t$, the same as those
of the polynomials $z_l^{d_l}+t\phi_l(z_l)$ where $\phi_l$ is an  arbitrary
polynomial in $z_l$ of degree less than $d_l$ defined near $z_l=0$.
In fact, by the linearity of $D\Psi$, it suffices to do this for
      $\phi_l(z_l)=z_l^k$ for each $0\leq k<d_l$.

Choose  smooth bump functions $\beta_l$ supported in disjoint balls
around the zeros  of $f$ with $\beta_l\equiv 1$ in a neighborhood of
$z_l=0$. For simplicity we fix $l$ and omit it from the notation.
We also fix local coordinates $\{v_j\}$ for $V$ around $f(0)$, and extend
these to coordinates $(v_j,x)$ for $X$ around $f(0)$, with $V$ given
locally as $x=x^1+Jx^2=0$ as described before Lemma \ref{ContactLemma}.

For any function $\eta(z)$  with $\eta(0)=1$, we can
construct maps
\bear
f_t=\l( \;\mspa f_0^T, \; \mspa f_0^N+t  \beta z^k \eta \;\r).
\label{5.ft}
\eear
It is easy to check that the zeros of the second factor have the form
$z_t(1+O(t))$ where the $z_t$ are the zeros of $z^d+tz^k$.
Then the variation $\dot f$ at time $t=0$ is $\xi= \beta z^k \eta \  e_N$,
where $e_N$ is a normal vector to $V$.

Keeping $x, p, j,J$ fixed, we will show that we can choose $\eta$ and a
variation $\dot \nu$ in $\nu$ such that $(0,\dot\nu,\xi,0, 0, 0)$ is
tangent to  $\UM_{g,n+\ell(s)}$.  This requires two conditions on 
$(\xi,\dot\nu)$.

\medbreak
(i)  The variation in $(J,\nu)$, which we are taking to  be $(0,\dot\nu)$,
must be tangent to ${\cal
J}^V$.  Thus $\dot\nu$ must satisfy the linearization of equations
(\ref{defJVeq}), namely
$$
\dot\nu^N=0 \qquad\mbox{and }\qquad    \left[\nabla_{e_N}\dot\nu+
J\nabla_{J e_N}\dot\nu\right]^N(\cdot )\ =\
\left[(J\,\nabla_{\dot\nu(\cdot )}J)e_N \right]^N
$$
along $V$, with $e_N$ as above.  This is true whenever
$\dot\nu$, in the coordinates of Lemma \ref{ContactLemma},  has an
expansion off $x=0$ of the form
\bear
\label{formofdotnu}
\dot\nu\ =\  A(z,v)+B(z,v)\,\bar{x}\,+\, O(|z|\,|x|)
\eear
with $A^N=0$ and $B^N=B^N(A)=\frac12[J(\nabla_{A(\partial/\partial z)}
J)(e_N)]^N$.

\medbreak

(ii)  If  $(0,\dot\nu,\xi,0, 0, 0)$ is to be tangent to the
universal 
moduli
space it must be in the
kernel of the linearized operator of Lemma  \ref{linearizationlemma}, and so must
satisfy
\bear
D\xi(z) - \dot \nu(z,f(z))=0
      \label{grapheqdotnu}
\eear
where $D$, which depends on $f$, is given in terms of the $\del$
operator of the pullback
connection by
$$
D\xi=\ov{\partial}_{f} \xi +\frac12(\nabla_\xi J) df\circ j-\nabla_\xi \nu.
$$

Near the origin in $(z,v,x)$ coordinates, (\ref{formofdotnu}) is a
condition on the 1-jet of
$\dot\nu^N$ along the set where $x=0$, and (\ref{grapheqdotnu}) is a
condition along the graph
$\{(z,v(z),z^d)\}$ of $f_0$.  Locally, these sets intersect only at the
origin.  Writing $\dot\nu=\dot\nu^V+\dot\nu^N$, we  take
$$
\dot\nu^V\ =\ [D\xi]^V
$$
along the graph and extend \pagebreak it arbitrarily to a neighborhood of the origin.
We can then take
$\dot\nu^N$ of the form (\ref{formofdotnu}) provided we can solve
\bear
D^N\xi(z) = \dot \nu^N(z, v(z),  z^d) = B^N(D^V\xi(0)) \,\bar{z}^d \,+\,
O(|z|^{d+1})
\label{dfexp1}
\eear
locally in a neighborhood of the origin with $D^N$ as in (\ref{3.4}).

Now, write $\xi=\alpha\,e_N+\beta \,Je_N$ where $\alpha$ and $\beta$ are
real, and identify
this with $\xi=\zeta\,e_N$ where $\zeta=\alpha+i\beta$ is complex.
Because $D$ is an $\R$-linear
first order operator, one finds that
\bear
D\xi & = &   (\del \zeta )\,e_N \ +\  \zeta E +\ov{\zeta}F
\label{firstDxi}
\eear
where $e_N$ is normal,
$$
E\ =\ \frac12\left[D(e_N)-J D(Je_N)\right]
\qquad\mbox{and}\qquad
F\ =\ \frac12\left[D(e_N)+J D(Je_N)\right].
$$
We need a solution of the form $\zeta=\beta z^k \eta$
near the origin.  For this we can take $\beta\equiv 1$.  The
equation (\ref{dfexp1}) we  must solve has the form
\begin{equation}
-z^k\del\eta \ =\ z^k\eta\,E^N(z,\ov z)\  + \bar{z}^k\bar\eta F^N(z,\bar z)
      + B^N([D\xi(0)]^V)\,\bar{z}^d \ +\ \ O\left(|z|^{d+1}\right).
\label{solvethisDxi}
\end{equation}
When $k=0$, (\ref{solvethisDxi}) has the form $\del\eta+a(z,\bar
z)\eta+b(z,\bar z)\bar\eta = G(z,\ov z)$, which
can always be solved by power series.  When $1\leq k<d$,  
we have $\zeta(0)=0$, so that
$B^N([D\xi(0)]^V)$  vanishes by
(\ref{firstDxi}).  Then  using Lemma
\ref{FNLemma} below, (\ref{solvethisDxi}) holds whenever $\eta$ satisfies
$$
-\del\eta \ =\ \eta\,E^N(z,\ov z)\  + a\bar{\eta}\bar{z}^kz^{d-1-k}
      \ +\ \ O\left(|z|^{d+1-k}\right),
$$
and this can also be  solved by power series. \enddemo

\proclaim{Lemma}
\label{FNLemma}
Near the origin, $F^N = a z^{d-1 }+O\left(|z|^{d}\right)$
for some constant~$a$.
\endproclaim

\pf Fix a vector $u$ tangent to the domain of $f$.  Using  the
definition of $F$, equation (\ref{linearizationequation}), and  the
$(J,\nu)$-holomorphic map equation $f_*u=2\nu(u)-Jf_*ju$, one finds that
$F^N(e_N)(u)= F^N(f_*u, u)$ where
\bear\qquad
4F^N(U, u)& = & J(\nabla_UJ)e_N - (\nabla_{JU}J) e_N +(\nabla_{e_N}J)JU\label{FNmess}
\\
&& +\
J(\nabla_{Je_N}J)JU + 2
(\nabla_{J\nu(U)}J)e_N -    2 (\nabla_{e_N}J)J\nu(u) \nonumber \\
&&-\ 2(\nabla_{Je_N}J)\nu(u) -2(\nabla_{e_N}\nu)u -
2J(\nabla_{Je_N}\nu) u. \nonumber
\eear
But the normal component of $U=f_*u$ is $dz^{d-1}\partial/\partial x$.  Thus
      we can replace $U$ in (\ref{FNmess}) by its  component in the $V$
direction;  the difference has
the form $z^{d-1} \Phi_1(z,\bar{z})$.  In the resulting expression, the $J$
is evaluated at the target
point: $J:=J(v(z),z^d)$.  But
$$
J(v(z),z^d)\ =\  J(v(z), 0) +O(|z|^d)
$$
and similarly  $\nabla J=(\nabla J)(v(z), 0)+O(|z|^d)$.  Finally, with $U$
tangent to $V$ and $J$
and $\nabla J$ replaced by their values at $(v(z), 0)$, one
can check that (\ref{FNmess}) vanishes by (\ref{defJVeq}). Lemma
\ref{FNLemma} follows.\enddemo

\section{Intersection data and rim tori}

The images of two $V$-regular maps can be distinguished by  (i)  their
intersection points with $V$, counted with
multiplicity, and (ii) their homology classes $A\in H_2(X)$.  One can go a
bit further:  if  $C_1$ and $C_2$ are
the images of two $V$-regular maps with the same data (i) and (ii), then
the difference $[C_1\#(-C_2)]$ represents
a class in $H_2(X\setminus V)$.  This section  describes  a space  ${\cal
H}^V_X$ of data that include   (i) and
(ii) plus enough additional data to make this last distinction.
Associating these data to a $V$-regular map then produces a continuous map
\best
\M^V_{g,n}(X) \to {\cal H}^V_X.
\eest
It is this map, rather than the simpler map to the data (i) and (ii), that
is needed for a gluing theorem for
relative invariants (\cite{IP4}).

We first need a space that records how $V$-regular maps intersect $V$.
Recall that the domain of each
$f\in\M^V_{g,n,s}$ has $n+\ell(s)$ marked points, the last
$\ell(s)$  of which are mapped into $V$.  Thus there is  an  {\it intersection map}
\bear
i_V: \M^V_{g,n,s}(X,A)\ra V_s
\label{5.intersectionmap}
\eear
that records the points and multiplicities where the image of $f$
intersects $V$, namely
\best
i_V\left(f,C, p_{1} \dots, p_{n+\ell}\right)\ =\
\left( \, (f(p_{n+1}),s_1),\dots, (f(p_{n+\ell}),s_\ell) \, \right).
\eest
Here $V_s$ is the space, diffeomorphic to  $V^{\ell(s)}$, of all sets of
pairs $ ( (v_1,s_1), \dots$,\break
$(v_\ell,s_\ell))$ with $v_i\in V$.  This is, of course, simply the
evaluation map at the last $\ell$ marked
points, but cast in a form that keeps track of multiplicities.

To simplify notation, it is convenient to take the union over  all
sequences $s$ to obtain the intersection map
\bear
i_V:    \M^V_{g,n}(X,A)\longrightarrow {\cal S}V
\label{5.bigintersectionmap}
\eear
where both
\bear
\M^V_{g,n}(X)\ =\ \coprod_{A}\coprod_{s}\  \M^V_{g,n,s}(X,A)
\qquad\mbox{and}\qquad
{\cal S}V\ =\
\coprod_{s}\ V_s
\label{4.defSV}
\eear
are given the topology of the disjoint union.

The next step is to augment ${\cal S}V$ with homology data to construct the
space ${\cal H}^V_X$. The  discussion
in the first paragraph of this section might suggest taking ${\cal H}$ to
be  $H_2(X\setminus V) \times {\cal
S}V.$  However,  the above images $C_1$ and $C_2$ do not lie in $X\setminus V$
--- only the difference does. In fact, the difference lies in
\bear
{\cal R}\ =\ {\cal R}_X^V\ =\ \mbox{ker}\ \left[ H_2(X\setminus V)\to
H_2(X)\right].
\label{defrimtori}
\eear
Furthermore,  there is a subtle twisting of these data, and  ${\cal H}$
turns out to be a nontrivial covering
space  over $H_2(X)\ti {\cal S}V$ with ${\cal R}$ acting as \pagebreak deck
transformations --- see (\ref{5.HSVproj})
below.  To clarify both these issues, we will  compactify
$X\setminus V$ and show how the images of $V$-regular maps determine cycles
in a homology theory for the
compactification.

Let $D(\ep)$ be the $\ep$-disk bundle in the normal bundle of $V$, identified with
a tubular neighborhood of $V$.
Choose a diffeomorphism of $X\setminus\overline{D(\ep)}$ with $X\setminus
V$ defined by the flow of a radial vector field and set
$S=\partial\overline{D(\ep)}$.  Then
\bear
\widehat{X}\ =\ \left[X\setminus \overline{D(\ep)}\right]\,\cup\,S
\label{hatX}
\eear
is a compact manifold with $\partial\widehat{X}=S$, and there is a
projection $\pi:\widehat{X}\to X$ which is the
projection $S\to V$ on the boundary and is a  diffeomorphism in the interior.

The appropriate homology theory is built from chains which, like the images
of $V$-regular maps, intersect $V$
at finitely many points.  Moreover, two cycles are homologous when they
intersect $V$ at the same points and
      their difference is trivial in  $H_2(X\setminus V)$.  We will  give two
equivalent descriptions of this
homology theory.

For the first description, consider

\begin{itemize}

\item[(i)] the free abelian group $C_k$ on
$k$-dimensional simplices in $\widehat{X}$, and

\item[(ii)] the subgroup $D_k$ generated by the $k$-simplices that lie entirely in one
circle fiber of $\partial\widehat{X}$.

\end{itemize}

\noindent Then  $({C_*}/{D_*}, \partial)$ is a chain complex over $\Z$.  Let
${\cal H}$ denote the 2-dimensional homology of this complex.  Elements
of $H_1(D_*)$ are linear
combinations of the circle fibers of $\partial\widehat{X}$.  Hence
$H_1(D_*)$ can be identified
with the space  ${\cal D}$   of divisors on $V$ (a  {\it divisor} is a finite
set of points in $V$, each with sign and multiplicity).  The long exact
sequence of the pair
$({C_*},{D_*})$  then becomes, in part,
\begin{equation}
0 \longrightarrow H_2(\widehat{X}) \stackrel{\iota}{\longrightarrow}
H_2({C_*}/{D_*})
\stackrel{\rho}{\longrightarrow}  {\cal D}.
\label{Hlongexact}
\end{equation}

For the second description we change the topology on $X$ and $\widehat{X}$
to separate cycles whose intersection
with $V$ is different.   Let $V^*$ be $V$ with the discrete
topology, and let
$S^*$  be
$S$ topologized as the disjoint
union of its fiber circles.  Then
$\pi:S^*\to V^*$ and the inclusions $V^*\subset X$ and
$S^*\subset \widehat{X}$ are continuous, and, when we use coefficients in $\Z$,
$H_1(S^*)$ is identified with the
space of divisors. The long exact sequence of the pair $(\widehat{X},S^*)$
again gives (\ref{Hlongexact}) with
$H_2(\widehat{X},S^*)$ in the middle.  To fix notation we will use this
second description.

The space in the middle of (\ref{Hlongexact}) is essentially the space of
data we want.  However, it is
convenient to modify it in two ways. First, observe that projecting
2-cycles into $X$ defines
maps $\pi_*:H_2(\widehat{X})\to H_2(X)$ and $\pi'_*:H_2(\widehat{X},S^*)\to
H_2(X)$ with $\pi_*=\pi'_*\circ
\iota$.  The kernel of $\pi_*$ is exactly the space ${\cal R}$ of
(\ref{defrimtori}), so that (\ref{Hlongexact}) can
be rearranged to read
\begin{equation}
0 \longrightarrow {\cal R} \stackrel{\iota}{\longrightarrow}
H_2(\widehat{X},S^*)
\stackrel{\rho}{\longrightarrow } H_2(X)\ti {\cal D}.
\label{Hlongexact2}
\end{equation}

Second, in keeping with what we have done with $V$-regular maps, we can
replace the space of divisors in
(\ref{Hlongexact}) by the space ${\cal S}V$ of (\ref{4.defSV}) which keeps
track of the numbering of the
intersection points, and whose topology  separates strata with different
multiplicity vectors $s$.  There is a continuous covering map ${\cal S}V\to {\cal D}$ which replaces ordered points by
unordered points. Pulling this covering back along the map
$\rho:H_2(\widehat{X},S^*)\to{\cal D}$ gives, at last, the desired space of
data.

\numbereddemo{Definition}
\label{5.defcalH}
Let ${\cal H}^V_X$ be the space $H_2(\widehat{X},S^*)
\ti_{\cal D} {\cal S}V$.
\enddemo

With this definition, (\ref{Hlongexact2}) lifts to a covering map
\begin{equation}\begin{array}{cccl}
{\cal R} & \longrightarrow &{\cal H}_X^V & \\
& & \Big\downarrow \ep&  \\
& &   H_2(X)\ti {\cal S}V. &
\end{array}
\label{5.HSVproj}
\end{equation}
where $H_2(X)$ has the discrete topology, ${\cal S}V$ is topologized as in
(\ref{4.defSV}), and
$\ep=(\pi'_*,\rho)$.  This is also the right space for keeping track of the
intersection-homology data:
      given a $V$-regular stable map $f$ in $X$ whose image is $C$, we
can restrict to $X\setminus V$, lift to $\widehat{X}$,
and take its closure, obtaining a curve $\widehat{C}$ representing a class $[\widehat{C}]$ in
$H_2(\widehat{X},S^*)$.  This is
consistent with the intersection map (\ref{5.bigintersectionmap}) because
$\rho[\widehat{C}]=\iota_V(f)\in{\cal D}$.
Thus there is a well-defined   map
\bear
h:\M_{g,n}^V(X)  \longrightarrow  {\cal H}_X^V
\label{5.defh}
\eear
which lifts the intersection map (\ref{5.bigintersectionmap})  through
(\ref{5.HSVproj}).  Of course, ${\cal H}^V_X$
has components labeled by $A$ and $s$, so this is a union of maps
\bear
h:\M_{g,n,s}^V(X,A)  \longrightarrow  {\cal H}_{X,A,s}^V
\label{5.defh2}
\eear
with $A\cdot V=\mbox{deg}\, s$.

We conclude with a geometric description of elements of ${\cal R}$ and of
the twisting in the covering
(\ref{5.HSVproj}).  Fix a small tubular neighborhood $N$ of $V$ in $X$ and
let $\pi$ be the projection from the
`rim' $\partial N$ to
$V$. For each simple closed curve $\gamma$ in $V$, $\pi^{-1}(\gamma)$ is a
torus in $\partial N$; such tori are
called {\it rim tori}.

\proclaim{Lemma}
 \label{rimtorilemma}
     Each element  $R\in {\cal R}$  can be represented by  a rim torus.
\endproclaim

\pf  Write $X$ as the union of $X\setminus V$ and a neighborhood  of $V$.
Then the
Mayer-Vietoris sequence
$$
\longrightarrow H_2(\partial{\hat X} )\stackrel{(\iota_*,\pi_*)}
\longrightarrow H_2(X\setminus V) \oplus H_2(V)
\longrightarrow H_2(X)  \longrightarrow
$$
shows that $(R,0)=\iota_*\tau$ for some $\tau\in H_2(\partial{\hat X} )$
with $\pi_*\tau=0$.  The lemma then
follows from  \pagebreak 
the  Gysin sequence for the oriented circle bundle \
$\pi:\partial{\hat X}\to V$:
\bear
\longrightarrow H_3(V)\stackrel{\psi}\longrightarrow
H_1(V) \stackrel{\Delta}\longrightarrow
H_2(\partial{\hat X} )\stackrel{\pi_*}\longrightarrow
H_2(V)\longrightarrow
\label{Gysin}
\eear
where $\psi$ is given by the cap product with the first  Chern class of the
normal bundle to $V$ in $X$.
\enddemo

Some rim tori are homologous to zero in $X\setminus V$ and hence do not
contribute to ${\cal R}$.  In
fact, the proof of Lemma \ref{rimtorilemma} shows that
$$
{\cal R}\ =\ \mbox{image}\ \left[\iota_* \circ \Delta:H_1(V) \to
H_2(X\setminus V)\right].
$$

Now consider the image $C$ of a $V$-regular map.  Suppose for simplicity
that $C$ intersects $V$ at a single
point $p$ with multiplicity one.  Choose a loop $\gamma(t)$, $0\leq t\leq
1$, in $V$ with $\gamma(0)=\gamma(1)=p$,
and let $R$ be the rim torus $\pi^{-1}(\gamma)$.  We can then modify $C$ by
removing  the annulus of radius
$\ep/2\leq r\leq \ep$ around $p$ in $C$ and gluing in the rim torus $R$,
tapered to have radius $\ep(1-t/2)$ over
$\gamma(t)$.   The resulting curve still intersects $V$ only at $p$, but
represents $[C]+[R]$.  Thus this gluing
acts as a deck transformation on $[C]\in {\cal H}_X^V$.  Retracting the
path $\gamma$, one also sees that
each ${\cal H}_{X,A,s}^V$ is path connected.

\numbereddemo{{R}emark}
\label{5.remark}
       There are no rim tori when $H_1(V)=0$ or when the map
$\iota_*\circ\Delta$ in (\ref{Gysin}) is zero.  In that
case
${\cal H}^V_X$ is simply $H_2(X)\ti {\cal S}V$.  In practice, this makes
the  relative invariants significantly
easier to deal with (see \S 9).
\enddemo

\section{Limits of $V$-regular maps}

In this and the next  section we  construct a compactification of each
component
       of the space of $V$-regular maps.  This compactification carries the
``relative virtual class'' that will enable us, in Section 8, to define the
relative GW
invariant.

One way to compactify $\M^V_{g,n,s}(X,A)$ is to take its closure
\bear
C\M^V_{g,n,s}(X,A)
\label{6.1}
\eear
in the space of stable maps $\ov\M_{g,n+\ell(s)}(X,A)$.
Under the `bubble convergence' of Theorem
\ref{BTCthm}  the limits of the last $\ell(s)$ marked points  are mapped
into $V$. Thus the closure lies in the
subset of $\ov\M_{g,n+\ell(s)}(X,A)$ consisting of stable maps whose last
$\ell(s)$ marked points are mapped into $V$; these still have
associated multiplicities $s_i$, although the actual order of contact
might be infinite.

The main step toward showing that
this closure carries a fundamental homology class is to prove that the
frontier $C\M^V\setminus \M^V$ is a
subset of codimension at least two.  For that, we  examine the elements of
$C\M^V$ and
characterize those stable maps that are limits of $V$-regular maps.  That
characterization allows us to count the
dimension of the
       frontier.  The frontier \pagebreak is a subset of the space of stable maps, so is
stratified
according to the type of bubble structure of the domain. Thus the goal
of this section is to work towards a proof of the
following statement about the structure of the closure $C{\M}^V$.

\proclaim{Proposition}  \label{compactificationThm}
For generic $(J,\nu)\in {\cal J}^V${\rm ,} each stratum of
the irreducible part of $$C\M^V_{g,n,s}(X,A)\setminus \M^V_{g,n,s}(X,A)$$
is an orbifold  of dimension at least two less than the dimension
{\rm (\ref{5.1})} of\break $\M^V_{g,n,s}(X,A)$.
\endproclaim

       The closure  $C\M^V$ contains strata corresponding to  different
types of limits.  For clarity these will be treated in several
separate steps:

\smallskip

Step 1: stable maps with no components or special points lying
entirely in $V$;

\smallskip

Step 2: a stable map with smooth domain which is mapped entirely  into~$V$;

\smallskip

Step 3: maps with some components in $V$ and some off $V$.

\demo{Step 1} For the strata consisting of stable maps with
no components  or special points in $V$ the analysis is essentially standard
(cf.\ \cite{RT1}).  Each stratum of this type is labeled by the genus
and the number $d\geq 1$ of double points of their nodal domain curve
$B$.  Fix such a $B$.  The corresponding stratum is the fiber of the
universal space $\pi:\UM^V_{B,n,s}(X,A) \to {\cal J}^V$ of\break $V$-regular
maps from $B$ into $X$, and  the irreducible part $\UM^{V*}$ of $\UM^V$ is
an orbifold by the same
tranversality arguments as in \cite{RT2}.
\enddemo

\proclaim{Lemma}
\label{7.Lemma2} In this {\rm `}\/Step {\rm 1'} case{\rm ,} for generic
$(J,\nu)\in {\cal J}^V${\rm ,} the irreducible part of the stratum
$\M_{B,n,s}^V(X,A)$ of $C{\M}^V$ is an  orbifold whose dimension
is $2d$ less than the dimension {\rm (\ref{Lemma5.1})} of $\M_{g,n,s}^V(X,A)$.
\endproclaim

\pf  Let $\wt B \ra B$  be the normalization of $B$.  Then $\wt B$ is
a (possibly disconnected) smooth curve with a pair of marked points
for each double point  of $B$.
We will show that $\UM^V_{B,n,s}(X,A)^*$  is a suborbifold of
$\UM^V_{g,n,s}(X,A)^*$ of codimension $2d$.   Lemma
\ref{7.Lemma2} then follows by the Sard-Smale theorem.

Assume for simplicity that there is only one pair of such marked points
$(z_1,z_2)$. Evaluation
at $z_1$ and $z_2$ gives a map
$$
\ev:\ \UM^V_{\wt B,n,s}(X,A)^* \ \to\  X\ti X
$$
and $\UM^V_{B,n,s}(X,A)^*$ is the inverse image of the diagonal $\Delta$ in
$X\ti X$.  Since
$\UM^V_{\wt B,n,s}(X,A)^*$ is an orbifold, we need only check that this
evaluation map is transversal to $\Delta$.

To that end, fix $(f_0,J,\nu)\in \ev^{-1}(\Delta)$.  Choose local
coordinates  in $X$ around $q=f_0(z_1)=f_0(z_2)$ and  cutoff functions
$\beta_1$ and $\beta_2$
supported in small disks around $z_1$ and $z_2$.  Then, as in (\ref{5.ft}),
we can  modify
$f_0$ locally around $z_1$ by $f_t=f_0+t\beta_1v$ and around $z_2$ by
$f_t=f_0-t\beta_2v$, and
modify $\nu$ to $\nu_t=\bar{\partial} f_t$ on the graph of $f_t$.  The
initial derivative of
this path is a tangent vector $w$ to $\UM_{\wt B}^V(X)$ with
$\mbox{ev}_*(w)=(v,-v)$.  Thus
$\ev$ is transversal to $\Delta$. \enddemo

 {\it Step} 2.  Consider  the strata
$C_1{\M}^V$ of $C{\M}^V$ consists of all maps with smooth domain
whose image is contained in $V$.  Such
maps lie in
$\M_{g,n+\ell(s)}(V,A)$, and it might seem that we can focus on
$V$ and forget about $X$.  But we are
only examining the subset
\best
C_1\M^V_{g,n,s}(X,A)\cap {\M}_{g,n+\ell(s)}(V,A) 
\eest
that lies in the closure of $\M_{g,n,s}^V(X,A)$.
The maps in this closure have a special
property, stated as  Lemma \ref{TaubesRenomlemma}.  This property
involves the linearized operator.

       For each $f\in \M_{g,n+\ell(s)}(V,A)$ denote by  $D^V$ the
linearization of the equation $\del f =\nu$ at the map $f$. Note that the
restriction map
\best
{\cal J}^V\ra {\cal J}(V)
\eest
that takes a compatible pair $(J,\nu)$  on $X$ to its restriction to $V$
is onto.  Then by Theorem \ref{5.1} for generic $(J,\nu)\in {\cal J}^V$ the
irreducible part of the moduli space
$\M_{g,n+\ell(s)}(V,A)$ is a smooth orbifold of (real) dimension
\bear
\ind\; D^V\ =\ -2K_V[A]+\left(\dim V-6\right)(1-g)+2n+2\ell(s).
\label{6.dimm(V)}
\eear

There are several related operators associated with the maps $f$ in
this moduli space.  First, there is the linearization $D_s^X$ of the
equation $\del f =\nu$; this acts on sections of $f^*TX$ that have
contact with $V$, described by the sequence $s$, and with index   given
by (\ref{Lemma5.1}).  Next, there is the operator $D^N$ obtained by
applying $D^X$ to vector fields normal to $V$  and then
projecting back onto the subspace of normal vector fields.
Completion in the Sobolev space with  $m$ derivatives in $L^2$
gives a bounded operator
\bear
D^N:L^{m,2}(f^*N_V)\ra L^{m-1,2}(T^*C\otimes f^* N_V)
\label{DN}
\eear
which is $J$-linear by Lemma \ref{DNcxlemma}.  For $m>\mbox{deg}\, s$
the sections that
satisfy the linearization of the  contact conditions specified by
$s$ form  a closed $J$-invariant subspace $L^{m,2}_s(f^*N_V)$.   Let
$D^N_s$  denote the restriction  of  $D^N$ to that subspace
$L^{m,2}_s$. The index of $D_s^N$ is the index of $D_s^X$ minus
the index of
$D^V$, so that
\bear
\ind \;D^N_s\ =\ \ 2(c_1(N_V)[A]+1-g-\deg s)\;=\;2(1-g)
\label{indexDNs}
\eear
since $\deg s=A\cdot V= c_1(N_V)[A]$.

\proclaim{Lemma}  \label{TaubesRenomlemma}
 Each element of the closure $C{\M}_s^V(X)$
whose image is a single component that lies entirely in $V$ is a map
with $\ker D_s^N\neq 0$.
\endproclaim

\pf This is seen by a renormalization argument similar to one in
\cite{T}. Suppose that  a sequence
$\{f_n\}$ in $\M^V_s(X)$ converges to  $f\in \M(V)$; in the present
case there is no bubbling, so that
$f_n\to f$ in
$C^\infty$. For large $n$, the images of the maps $f_n$
lie in a neighborhood of $V$, which we identify with a
subset in the normal bundle $N_V$ of $V$ by the exponential map. 
Let $\phi_n$ be the
projection of $f_n$ to $V$ along the fibers of $N_V$, so that $\phi_n\to
f$ in $C^\infty$.

Next let $R_t:N_V\ra N_V$ denote the dilation by a
factor of $1/t$.  Because the image  of $f_n$ is not contained in $V$
there is, for each $n$,  a
unique $t=t_n$ for which the normal component of the pullback map
$R_t(f_n)$ has $C^1$ norm equal to
1.  These $t_n$ are positive and $t_n\to 0$.   Write $R_{t_n}$ as
$R_n$ and consider the  renormalized
maps $F_n=R_{n}(f_n)$. These are holomorphic with respect to renormalized
$(R_n^* J, R_n^* \nu)$; that is,
\bear
\del_{j,R_n^* J} F_n- R_n^* \nu=R^*_n( \del_{j,J}f_n-\nu)=0.
\label{6.3a}
\eear
By expanding in Taylor series one sees that $(R_n^* J,R_n^* \nu)$
converges in $C^\infty$ to a limit $(J_0,\nu_0)$; this limit is
dilation invariant and equal to the restriction of
$(J,\nu)$ along $V$.  The sequence $\{F_n\}$ is also bounded in
$C^1$.  Therefore, after applying  elliptic bootstrapping and passing
to a subsequence, $F_n$ converges
in $C^\infty$ to a limit $F_0$ which satisfies
$$
\del_{j,J_0} F_0-\nu_0=0.
$$
We can also write $F_n$ as $\exp_{\phi_n}\xi_n$ where
$\xi_n\in\Gamma(\phi_n^*N_V)$ is the normal component of $F_n$, which
has $C^1$ norm equal to 1.  The above convergence implies that $\xi_n$
converges in $C^\infty$ to some nonzero $\xi\in\Gamma(f_0^*N_V)$. We
claim that $\xi$ is in the kernel of $D^N_s$ along $f_0$.  In fact,
since the $f_n$ satisfy the contact constraints described by $s$ and
converge in $C^\infty$ the limit $\xi$ will have zeros described by
$s$. Hence we need only show that $D^N_{f_0}\xi=0$.

For  fixed $n$, $\phi_n$ and $f_n=\exp_{\phi_n}(t_n\xi_n)$ are maps from
the same domain so that by the definition of the linearization (for fixed
$J$ and $\nu$)
$$
P_n^{-1}(\del_J f_n-\nu_{f_n}) -(\del_J \phi_n-\nu_{\phi_n}) \ =\
D_{\phi_n}(t_n\xi_n,0) +
O\left(|t_n\xi_n|^2\right)
$$
where $P_n$ is the parallel transport along the curves
$\exp_{\phi_n}(t\xi_n)$, $0\leq t\leq t_n$.  The first term in this
equation vanishes because $f_n$ is
$(J,\nu)$-holomorphic.  Furthermore,   because the image of
$\phi_n$ lies in $V$, condition (\ref{defJVeq}a) means that the normal
component of
$\del_J \phi_n-\nu_{\phi_n}$ vanishes.  After dividing through by
$t_n$ and noting that
$t_n^{-1}|t_n\xi_n|^2 \leq t_n$ we obtain
\vglue12pt
\hfill ${\displaystyle
D^N_{f_0}\xi\ = \ \lim_{n\ra \infty} D^N_{\phi_n}(\xi_n,0)\ =\ 0.
}$\enddemo

The operator $D^N_s$ depends only on the 1-jet of $(J,\nu)\in {\cal J}^V$,
so that we can consider the restriction map
\bear
{\cal J}^V\ra {\cal J}^1
\label{6.Jrestriction}
\eear
that takes a compatible pair $(J,\nu)$  on $X$ to its 1-jet along $V$.
This map is onto, and by Lemma \ref{DNcxlemma},   $D^N_s$ is a
complex operator for any $(J,\nu)\in {\cal J}^1$. Then  $D^N_s$  defines a
smooth section of
\bear
\begin{array}{c}
{\rm Fred} \\[5pt]
\big\downarrow\\[5pt]
\UM(V)
\end{array}
\label{diagramFred}
\eear
where $\UM(V)$ is the universal moduli space of
maps  into $V$ (which is a fiber bundle over ${\cal J}^1$), and where
${\rm Fred}$ is the bundle whose fiber at $(f,j,J,\nu)$ is the space of
all  complex linear Fredholm  maps (\ref{DN}) of index $\iota\le 0$.
By a theorem of Koschorke \cite{K}, ${\rm Fred}$ is a disjoint union
$$
{\rm Fred} = \bigcup_k\, {\rm Fred}_k\
$$
where ${\rm Fred}_k$ is the complex codimension $k(k-\iota)$ submanifold
consisting of all  the operators whose kernel is exactly $k$ complex
dimensional. In fact,  the normal bundle to ${\rm Fred}_k$ in
${\rm Fred}$, at an operator $D$, is $\mbox{Hom}(\mbox{ker }
D,\mbox{coker }D)$.

\proclaim{Lemma}
\label{transverseLemma1}
The section $D_s^N$ of {\rm (\ref{diagramFred})} is transverse to
each ${\rm Fred}_k$.
\endproclaim
\pf  Fix $(f,j,J,\nu)\in \UM_s^V(X)$ such that the linearization
$D_s^N$ at $(f,j,J,\nu)$ lies on
${\rm Fred}_k$.  Let $\pi^N$ be the projection onto the normal part, so that
$D_s^N=\pi^N\circ D_s^X$.  The lemma follows if we show that for any
elements $\kappa\in \mbox{ker } D_s^N$ and
$c\in \mbox{ker } \l(D_s^N\r)^*$ we can find a variation in $(J,\nu)$
       such that
$$
\left\langle c,\  (\delta D_s^N) \kappa\right\rangle \neq 0
$$
(these brackets mean the $L^2$ inner product on the domain $C$ and
$\l(D_s^N\r)^*$ is the $L^2$ adjoint of
$D_s^N$). But
$$
(\delta D_s^N) \kappa \ =\  (\delta
\pi^N)  D_s^X \kappa +
\pi^N (\delta D^X) \kappa+ \pi^N D^X(\delta \pi^N) \kappa
$$
with the linearization $D^X$ is given by
(\ref{linearizationequation}). We will take the variation with
$(f,j,J)$ fixed and $\nu$ varying as $\nu_t=\nu+t\mu$ with $\mu\equiv
0$ along $V$. Then $\pi^N$ is fixed; i.e., it depends on $J$ and $f$, but
not on $\nu$.  Hence the above reduces to
\begin{equation}
\left\langle c, \ (\delta D^X) \kappa\right\rangle\ =\   -
\left\langle c,\;\nabla_\kappa \mu
\right\rangle.
\label{variationTian}
\end{equation}
This depends only on the 1-jet in the second variable of $\mu$ along $V$,
where $\mu$ is the variation in $\nu(x,f(x))$.

Choose a point $x\in C$ such that $\kappa(x)\neq 0$.  Let $W$ be a
neighborhood of $x$ in $\P^N$ and $U$ a
neighborhood of $f(x)$ in $X$  such that $\kappa$ has no zeros in $U$.  To
begin, $c$ is defined only
along the graph of $f$ and is a $(0,1)$ form with values in $N_V$.  Extend
$c$ to a smooth section
$\tilde{c}$ of $\mbox{Hom}(T\P^N,TX)$ along $W\ti U$ such that
$\tilde{c}|_V$ is a section of
$\mbox{Hom}(T\P^N,N_V)$.  Multiply $\tilde{c}|_V$ by a  smooth bump
function $\beta$ supported on $W\ti U$ with
$\beta\equiv 1$ on a slightly smaller open set.

Now construct the $(0,1)$ form $\mu$ such that its  1-jet  along $V$ 
satisfies
$$
\mu|_V=0, \enspace   
\nabla_{\kappa(y)} \mu(x,y) \, = \, \left(\beta\tilde{c}\right)(x,y)
\enspace\mbox{ and }\enspace
\nabla_{J\kappa(y)} \mu(x,y) \, =\,  - J\left(\beta\tilde{c}\right)(x,y).
$$
The required compatibility conditions
(\ref{defJVeq}) are now satisfied
because the right-hand side of (\ref{defJVeq}c) vanishes since $\mu$ vanishes
along $V$.  Moreover,
$$
\left\langle c, \ (\delta D_s^X) \kappa\right\rangle
\ =\ - \int_C \left\langle c, \ \nabla_\kappa \mu\right\rangle
\ =\ - \int_{C\cap U} \beta\,|c|^2.
$$
But $c$ satisfies the elliptic equation $(D_s^N)^*c=0$, so by the unique
continuation theorem for elliptic operators $|c|$ does not identically
vanish on any
open set.  Thus we
have found a nonzero variation. \enddemo

\proclaim{Proposition} \label{M'codim2}
In this {\rm `}\/Step {\rm 2'} case  $C_1\M^V_{g,n,s}(X,A)$ is
contained in the  space
\bear
\M'_{g,n,s}=
\left\{ (f,j)\in \M_{g,n+\ell(s)}(V,A)\ |\ {\mbox{\rm dim}} \
{\mbox{\rm ker}} \ D_s^N\ne 0\
       \right\}.
\label{6pnewspace}
\eear
Moreover{\rm ,} for generic $(J,\nu)\in {\cal J}^V${\rm ,} the irreducible part
of {\rm (\ref{6pnewspace})} is a suborbifold
of $\M_{g,n+\ell(s)}(V,A)^*$ of dimension two  less
than {\rm (\ref{5.1}).}
\endproclaim 

\pf The first statement follows from Lemma \ref{TaubesRenomlemma}.
Next, note that the dimension (\ref{6.dimm(V)}) of
$\M_{g,n+\ell(s)}(V,A)$ differs from
(\ref{5.1}) by exactly  the
index (\ref{indexDNs}) of $D^N_s$, so the second statement is
trivially true if
$\mbox{index } D^N_s >0$.
Thus we assume that $\iota=\mbox{index } D^N_s \leq 0$.

       Lemma \ref{transverseLemma1} implies that
the set of pairs $(f,j,J,\nu)\in \UM_{g,n+\ell(s)}(V)$ for which
$D_s^N$ has a nontrivial kernel, namely
$$
\UM' = D^{-1}\left({\rm Fred}\setminus{\rm Fred}_0\right),
$$
is  a  (real) codimension $2(1-\iota)$ subset of
$\UM_{g,n+\ell(s)}(V)$, and in fact  a suborbifold  off a set of
codimension $4(2-\iota)$.  Since  the projection
$\pi: \UM'  \to {\cal J}^1$ is Fredholm, the Sard-Smale theorem implies that
       for a  second category set of $J\in {\cal J}^1$ the fiber
$\pi^{-1}(J)$ --- which is the space
(\ref{6pnewspace}) --- is  an orbifold of (real) dimension
$$
2\,\mbox{index } D^V -2(1-\iota) \, =\,  2\,\mbox{index } D^V
+2\,\mbox{index } D^N_s -2\, =\,
       2\,\mbox{index } D^X_s-2.
$$
The inverse image of this  second category set under (\ref{6.Jrestriction})
is a   second category set in ${\cal J}^V$.  Hence (\ref{6pnewspace})
is an orbifold for
generic $(J,\nu)\in {\cal
J}^V$, and has codimension
at least two in $C\M^V_s(X)$. \enddemo

{\it Step} 3.  Next consider limit maps $f\in C{\M}_s^V(X)$
whose domain is the union $C=C_1\cup C_2$ of bubble domains of genus
$g_1$ and $g_2$ with $f$ restricting to a $V$-regular map $f_1:C_1\ra
X$ and a map $f_2:C_2\ra V$ into $V$. Limit maps $f$ of this type
arise, in particular, from sequences of maps in which either (a) two
contact points collide in the domain or (b) one of the original $n$
marked points collides with a contact point because its image sinks
into $V$.  In either case the collision produces a ghost bubble map
$f_2:C_2\to V$ which has energy at least $\al_V$ by Lemma
\ref{minenergylemma}.

\vglue12pt

In this Step 3 case, $f_1^{-1}(V)$ consists of the nodal points
$C_1\cap C_2$ and some of the last $\ell(s)$ marked points $p_k\in C$.
The nodes are defined by identifying points $x_j\in C_1$ with $y_j\in
C_2$.  Since $f_1$ is $V$-regular and $f_1(x_j)\in V$, Lemma
\ref{ContactLemma} associates a
multiplicity $s'_j$ to each $x_j$. Similarly, since $f$ arises as a
limit of $V$-regular maps the
$p_i$, being limits of the contact points with $V$, have associated
multiplicities.  The set of $p_i$ is split into the points $\{p_i^1\}$
on $C_1$ and $\{p_i^2\}$ on $C_2$; let $s^1=(s^1_1, s^1_2,\dots )$ and
$s^2=(s^2_1, s^2_2,\dots )$ be the associated multiplicity vectors.
Thus $f$ is a pair
\bear \qquad\quad
f= (f_1,f_2)\,\in\,\M^V_{g_1,n_1,s^1\cup
s'}(X,[f_1])\ti \ov{\M}_{g_2,n_2+\ell(s^2)+\ell(s')}(V,[f_2])
\label{6.step3.1}
\eear
with $n_1+n_2=n$, $[f_1]+[f_2]=A$, $\deg s^1+\deg s'=[f_1]\cdot V$,
and satisfying the matching
conditions $f_1(x_j)=f_2(y_j)$.

\proclaim{Proposition} \label{TaubesRenomlemma2}
   In this {\rm `}Step {\rm 3'} case{\rm ,} the only elements
{\rm (\ref{6.step3.1})} that lie in  $C{\M}_{s}^V(X)$ are
those for which there is a {\rm (}\/singular\/{\rm )} section
$\xi\in\Gamma(f_2^*N_V)$ nontrivial on at least one component of
$C_2$ with  zeros of order
$s^2_i$ at $p_i^2${\rm ,}  poles of order $s_j'$ at
$y_j$ {\rm (}\/and nowhere else\/{\rm ),} and $D_{f_2}^N\xi=0$ where $D_{f_2}^N$ is as in
{\rm (\ref{3.4}).}
\endproclaim  

The proof uses a renormalization argument similar to the one used in Lemma
\ref{TaubesRenomlemma}, but this time done in a compactification
$\P_V$ of  the normal bundle $\pi:  N_V\to V$.  For clarity we
describe $\P_V$ before
starting the proof.

      Recall that $N_V$ is a complex
line bundle with an inner product and a compatible connection induced
by the Riemannian connection on $X$.
      As a manifold $\P_V$ is the
fiberwise complex projectivization of the Whitney sum of $N_V$ with
the trivial complex line bundle
$$
{\pi_\P}:\P_V=\P(N_V\oplus\cx)\ \to\ V.
$$
Note that the bundle map  $\iota:N_V\hookrightarrow \P_V$ defined by
$\iota(x)=[x,1]$ on each fiber  is   an embedding
      onto the complement of the infinity section $V_\infty\subset\P_V$.
The scalar multiplication map $R_t(\eta)=\eta/t$ on
$N_V$ defines a $\cx^*$ action on $\P_V$.

      When $V$ is a point we can identify $\P_V$ with $\P^1$ and give it
the K\"{a}hler structure $(\w_\ep,g_\ep,j)$ of the
2-sphere of radius $\ep$.  Then
$\iota:\cx\to \P_V$ is a holomorphic map with
$\iota^*g_\ep = \phi_\ep^2\,\l[(dr)^2+r^2(d\theta)^2\r]$ and
$\iota^*\w_\ep = \phi_\ep^2\,r dr\wedge d\theta = d\psi_\ep\wedge
d\theta$ where
$$
\phi_\ep(r)  = \frac{2\ep}{1+r^2}\qquad \mbox{ and }\qquad
\psi_\ep(r) = \frac{2\ep^2r^2}{1+r^2}.
$$
This construction globalizes by interpreting $r$ as the norm on the
fibers of $N_V$, replacing $d\theta$ by the connection 1-form
$\alpha$ on $N_V$ and including the curvature  $F_\alpha$ of that
connection.  Thus
$$
\iota^*\w_\ep\ =\ \pi^*\w_V + \psi_\ep \pi^* F_\alpha + d\psi_\ep\wedge \al
$$
is a closed form which is nondegenerate for small $\ep$ and whose
restriction to each fiber of $N_V$  agrees with the volume form
on the 2-sphere of radius $\ep$.   Furthermore, at
each point $p\in N_V$  the connection determines a horizontal
subspace which identifies $T_pN_V$ with
the fiber of $N_V\oplus TV$ at $\pi(p)$.  But the fibers of $N_V$
have a  complex structure $j_0$ and a metric $g_0$, and
$J$ and $g$ on $X$ restrict to $V$. One can then check that for small
$\ep$ the form $\w_\ep$,
$$
\tilde{J} = j_0 \oplus J|_V,
\qquad\mbox{and}\qquad
\tilde{g}_\ep = \l(\phi_\ep^2 g_0\r) \oplus g|_V
$$
extend over $V_\infty$  to define a tamed triple
$(\w_\ep,\tilde{J},\tilde{g}_\ep)$  on $\P_V$. As in \cite{PW}, Lemma 
\ref{minenergylemma} holds for tamed
structures,  and so we can choose $\ep$ small enough  that every
$J$-holomorphic map
$f$ from $S^2$  onto a  fiber of $\P_V\to
V$ of degree $d\leq [f_1]\cdot V$ satisfies
\bear
\int_{S^2}  |df|^2\ \leq\ \frac{\alpha_V}{8}
\label{6.smallfiber}
\eear
where   $\alpha_V$ is the  constant associated with $V$  by Lemma
\ref{minenergylemma}. We fix
such an $\ep$ and write $\w_\ep$ as $\w_\P$. Let $V_0$ denote the zero section of $\Bbb P_V$.

Now symplectically identify an $\ep'<\ep$ tubular neighborhood of
$V_0$ in $\P_V$ with a   neighborhood of $V\subset X$ and
pullback $(J,g)$ from $X$ to $\P_V$. Fix a bump function $\beta$
supported on the $\ep'$ neighborhood of $V_0$ with $\beta=1$ on the
$\ep'/2$ neighborhood. For each small $t>0$ set $\beta_t=\beta\circ
R_{{t}}$.  Starting with the ``background'' metric $g'=
\beta_t g +(1-\beta_t) g_\ep$,  the procedure described in the
appendix produces a  compatible triple $(\w_\P,J_t, g_t)$ on $\P_V$.
Then as $t\to 0$ we have $J_t\to\tilde{J}$ in $C^0$ on $\P_V$  and
$  g_t \to g_0$ on compact sets of $\P_V\setminus
V_\infty$.

\demo{Proof of Proposition {\rm \ref{TaubesRenomlemma2}}}
Suppose that a sequence of $V$-regular maps\break $f_m:C_m\to X$ converges
to $f=(f_1, f_2)$ as above.  That means that the domains $C_m$
converge to $C=C_1\cup C_2$ and, as in Theorem \ref{BTCthm}, the $f_m$
converge to $f:C\to X$ in $C^0$ and in energy, and $C^\infty$ away
from the nodes of $C$.

Around each node $x_j=y_j$ of $C_1\cap C_2$ we have coordinates $(z_j,
w_j)$ in which $C_m$ is locally the locus of $z_jw_j=\mu_{j,m}$ and
$C_1$ is $\{z_j=0\}$.  Let $A_{j,m}$ be the annuli in the neck of
$C_m$ defined by $|\mu_{j,m}|/\delta\leq |z_j|\leq \delta$.  We also
let $C_m'\subset C_m$ denote the neck $A_m=\cup_j A_{j,m}$ together
with everything on the $C_2$ side of $A_m$, $f'_m$ be the restriction
of $f_m$ to $C'_m$, and let $\phi_m$ be the corresponding map into the
universal curve as in (\ref{2.defjholo}).

The restrictions of $f_m$ to  $C'_m\setminus A_m$ converge to $f_2$.
Because the image of $f_2$ lies in $V$ its energy  is at
least the constant $\alpha_V$ associated with $V$  by Lemma~\ref{minenergylemma}.  We can then fix $\delta$ small enough
that the energy of $f=(f_1, f_2)$  inside the union of $\delta$-balls around
the nodes is at most $\alpha_V/32$.  Then  for  \pagebreak large~$m$
\begin{eqnarray}
 \int_{C'_m\setminus A_m} |d(\pi_V\circ  f_m)|^2 + |d\phi_m|^2& \geq& \alpha_V/2
 \label{6.6.1}\\
\noalign{\noindent and}
\int_{A_m} |df_m|^2+ |d\phi_m|^2&\leq& \alpha_V/16.
\nonumber
\end{eqnarray}

To renormalize, note that for large $m$ the image  of  $f_m'$ lies in
a  tubular neighborhood of $V$ which is identified with a
neighborhood of $V_0$ in $\P_V$.   Hence  $f'_m$ gives rise to a
one-parameter family of maps $\iota\circ R_t\circ
f_m'$ into $\P_V$.   We can consider the energy (\ref{defEpfp}) of
the corresponding map $(\iota\circ R_t\circ f_m',\,\phi_m):
C_m'\to \P_V\times\ov{\cal U}_{g,m}$ on the part of the domain which
is mapped into the upper
hemisphere $\P_V^+$ calculated using the metric $\wt g_\ep$ on $\P_V$
constructed above. That energy vanishes for large $t$ and exceeds
$\alpha_V/4$ for small
$t$ by (\ref{6.6.1}).  Therefore there is a unique
$t=t_m$ such that the maps
$$
g_m:C_m'\to\P_V \qquad\mbox{by}\qquad g_m= \iota\circ R_{t_m}\circ f_m 
$$
satisfy
\bear
\int_{g_m^{-1}\l(\P^+_V\r)\cup A_m} |dg_m|^2 + |d\phi_m|^2\ =\  \al_V/4.
\label{6.6albound}
\eear
Note that $t_m\to 0$ because of (\ref{6.6.1}) and the fact that
      $f_m(C_m'\setminus A_m)\to V$ pointwise.

Next consider the small annuli $B_{j,m}$ near $\partial C'_m$ defined
by $\delta/2\leq |w_j|\leq \delta$ and let $B_m=\cup_jB_{j,m}$. On
each $B_{j,m}$ $f_m$ converges in $C^1$ to $f_1=a_jw_j^{s_j}+\dots $
and $f_m(B_{j,m})$ has small diameter.  Hence, after possibly making
$\delta$ smaller and passing to a subsequence, each $g_m(B_{j,m})$
lies in a coordinate neighborhood $V_j$ centered at a point $q_j\in
V_\infty$ with $\mbox{diam}^2(V_j)<\alpha_V/1000$.  Fix a smooth bump
function $\beta$ on $C_m$ which is supported on $C_m'$, satisfies $
0\leq \beta\leq 1$ and $\beta\equiv 1$ on $C'_m\setminus B_{m}$, and
so that the integral of $|d\beta|^2$ over each $B_{j,m}$ is
bounded by $100$.

Now extend $C_m'$ to a closed curve by smoothly attaching a disk $D_j$
along the circle $\gamma_{j,m}=\{|w_j|=\delta\}$.  Extend $g_m$ to
$\ov{g}_m:\ov{C}_m=C_m'\cup \{D_j\}\to\P_V$ by setting
$\ov{g}_m(D_j)=q_j$ and coning off $g_m$ on $B_{j,m}$ by the formula
$\ov g_m=\beta\cdot g_m$ in the coordinates on $V_j$. The local
expansion of $f_1$ shows that $f_m(\gamma_{j,m})$, oriented by the
coordinate $w_j$, has winding number $s_j$ around $V_0$.  The same is
true of $g_m(\gamma_{j,m})$, so in homology $[\ov{g}_m]$ is
$\iota_*[f_2]+sF$ where $s=\sum s_j$ and $F$ is the fiber class of
$\P_V\to V$.

By (\ref{6.6albound}), the  energy of $\ov{g}_m$ on the region
that is mapped into $\P_V^+$ is  bounded by
\bear
\label{capenergybound}
\int_{g_m^{-1}\l(\P^+_V\r)} \l|dg_m\r|^2 \ +\ \sum_j
\mbox{diam}^2(V_j) \,\int_{B_{j,m}} |d\beta|^2
\ \leq\ \frac{\alpha_V}{2}.
\eear
On the other hand, in the region  mapped into $\P_V^-$,
$\ \ov{g}_m=g_m$ is $(J_m,\nu_m)$-holomorphic with $J_m\to
\tilde{J}$ and $\nu_m\to\pi^*\nu_V$,
so the energy in  that region   is dominated by its symplectic
area (1.6).  Thus
\begin{eqnarray*}
E(\ov{g}_m)
&\leq& \frac{\alpha_V}{2}  \ +\ c_1\,\int_{g_m^{-1}\l(\P^-_V\r)}g_m^*\w_\P
\\ &\leq&     \frac{\alpha_V}{2} \ +\ c_1\,\l\langle \w_\P,\
[\ov{g}_m]\r\rangle \ +\ c_1\,\int_{g_m^{-1}\l(\P^+_V\r)}
\l|\ov{g}_m^*\w_\P\
\r|.
\end{eqnarray*}
With (\ref{capenergybound}) this gives a uniform energy bound of the
form $E(\ov{g}_m)\leq   c_1\,\langle\w,\break
[f_2]+sF\rangle+c_2 $.

      This energy bound applies, {\it a fortiori}, to the restrictions
$g_m'$ of $g_m$ to $C'_m\setminus B_m$.  These $g_m'$ are
$(J_m,\nu_m)$-holomorphic, so Theorem \ref{BTCthm} provides a subsequence
which converges to a $(\tilde{J},\pi^*\nu_V)$-holomorphic map whose domain
is
$C_2$ together with the disks $\{|w_j|\leq \delta/2\}$ in $C_1$ and
possibly some bubble components.

After deleting those disks, the limit is a map
$g_0:\tilde{C}_2\ra \P_V$ with $g_0(y_j)\break\in
V_\infty$ at marked points $y_j$.   By construction, the projections
$\pi\circ g_m'$ converge to $f_2$, so the irreducible components of
$\tilde{C}_2$ are of two types: (i) those  biholomorphically identified
with  components of $C_2$ on which
$g_0$ is a lift of $f_2$ to
$\P_V$, and (ii) those
mapped by $g_0$ into fibers of $\P_V$ and also collapsed by the
stabilization
$\tilde{C}_2\mapsto {\rm st}(\tilde{C}_2)$   Then (\ref {6.6albound})
implies that no
type (i) component is mapped to
$V_\infty$.  The type (ii) components are $(J,0)$-holomorphic and on them
$|d\phi|^2\equiv 0$, so by (\ref{6.smallfiber}) these components
contribute a total of at most $\alpha_V/8$ to the integral (\ref
{6.6albound}).  Thus (\ref{6.smallfiber}) implies that at least one
component of type (i) is not mapped into $V_0$.

       Lemma \ref{ContactLemma} shows each component of $g_0$ has a
local expansion normal to $V_\infty$ given by $b_j z_j^{d_j}+\cdots \
$ at each $y_j$.  To identify $d_j$ we note that $\partial
A_{j,m}=\gamma_{j,m} \cup \gamma_{j,m}'$ where $\gamma_{j,m}'$ is the
circle $|z_j|=\delta$ oriented by $z_j$.  The homology
$g_m(A_{j,m})\subset \P_V\setminus V_\infty$ then shows that $d_j$,
which is the local winding number of $g_m(\gamma'_{j,m})$ with
$V_\infty$, is equal to the local winding number of
$g_m(\gamma_{j,m})$ with $V_\infty$, which is $s_j$.

The convergence $g'_m\to g_0$ on $C_2$ means that the sections
$\xi_m=\iota^{-1}g_m$ of $f_2^*N_V$ converge to a nonzero
$\xi=\iota^{-1}g_0$.  Then $D^N\xi=0$ as in the proof of Lemma
\ref{TaubesRenomlemma}, and our intersection number calculation shows
that $\xi$ has a pole of order $s_j$ at each node $y_j$.  Furthermore,
the $g_m$ have the same zeros, with multiplicity, as the $f_m$, so the
zeros of $\xi$ are exactly the last $\ell(s)$ marked points of the
limit curve $C_2$ and the multiplicity vector associated with those
zeros is the original $s$.  Thus $\xi$ is a nonzero element of $\ker
D_{s,s'}^N$.  \enddemo

Proposition \ref{TaubesRenomlemma2} shows   that maps of the form
(\ref{6.step3.1}) which are in the closure of $C{\M}_{s}^V(X)$ carry a
special structure: a nonzero element $\xi$ in the kernel of $D^N_{f}$
with specified poles and zeros, defined on some component that is
mapped into $V$.  That adds constraints which enter the dimension
counts needed to prove Proposition \ref{compactificationThm}.  In fact
the proof shows that $\xi$ vanishes only on those components which
sink into $V_0$ as the renormalized maps $g_n$ converge.  On those
components we can renormalize again and proceed inductively.  But
instead of continuing down this road of special cases, we will define
the special structure in the general case of maps with many
components.  Those maps form a space of `$V$-stable maps', and we will
then do the dimension count once and for all in that context.

\section{The space of $V$-stable maps}

            In the general case, the limit of a sequence of $V$-regular
maps is a stable map whose components are of the types described in
Steps 1--3 of Section~6. The components of the limit map are also
partially ordered according to the rate at which they sink into $V$.
In this section we introduce terminology which makes this precise, and
then construct a compactification for the space of $V$-regular maps.

Let $C$ be a stable curve.  A {\it layer structure} on $C$ is the
assignment of an integer $\lambda_j=0,1,\dots$ to each irreducible
component $C_j$ of $C$.  At least one component must have
$\lambda_j=0$ or $1$.  The union of all the components with
$\lambda_j=k$ is the  {\it layer $k$ stable  curve} $B_k\subset
C$. Note that $B_k$ might not be a connected curve.

\numbereddemo{Definition}
A {\it marked layer structure} on $C\in\ov{\M}_{g,n+\ell}$ is  a layer
structure on $C$ together with
\begin{itemize} \item[(i)] a vector $s$ giving the multiplicities of the last
$\ell=\ell(s)$  marked
points, and

\item[(ii)] a  vector $t$ that assigns multiplicities to each double point of
$B_k\cap B_{l}$, $k\neq l$.
\label{7.markedlayerstructure}
\end{itemize}
\enddemo

 Each layer $B_k$ then has points $p_{k,i}$ of type (i) with
multiplicity vector $s_k=(s_{k,i})$, and has double points with
multiplicities.  The double points separate into two types.  We let
$t^+_k$ be the vector derived from $t$ that gives the multiplicities
of the double points $y_{k,i}^+$ where $B_k$ meets the higher layers,
i.e. the points $B_k\cap C_j$ with $\lambda_j>k$.  Let $t^-_k$ be the
similar vector of multiplicities of the double points $y_{k,i}^-$
where $B_k$ meets the lower layers.  Note that the double points
within a layer are not assigned a multiplicity.

There are operators $D_k^N$ akin to (\ref{DN}) defined on the layers
$B_k$, $k\geq 1$, as follows.  The marked points
$y^-_{k,i}$ define $\ell(t^-_k)$ disjoint sections of the   universal curve
$\ov\U_{g,n+\ell}\to \ov\M_{g,n+\ell}$;  in fact by
compactness those sections have  disjoint tubular neighborhoods. For
each choice of $t=t^-_k$ and $\alpha$, fix smooth  weighting
functions $W_{t,\alpha}$ whose restriction to each fiber of the
universal curve has the form $|z_j|^{\alpha+t^-_{k,j}}$ in
some local coordinates
$z_j$ centered on $y_{k,j}^-$ and has no other zeros.  Then given a
stable map  $f:B_k\to V$  let $L^m_{t,\delta}(f^*N_V)$ be the
Hilbert space of all
$L^m_{\rm loc}$ sections $f^*N_V$ over
$B_k\setminus
\{y^-_{k,j}\}$ which are finite in the norm
$$
\|\xi\|^2_{m,t,\delta}\ =\ \sum_{l=0}^m\int_{B_k}\left|
W_{t,l+\de}\cdot\nabla^l \xi\right|^2.
$$
For large $m$ the elements $\xi$ in this space have poles with
$|\xi|\leq c |z_j|^{-t^-_j-\delta}$ at each $y^-_{k,j}$ and have $m-1$
continuous derivatives elsewhere on $B_k$.
For such $m$  let  $L^m_{k,\delta}(f^*N_V)$ be  the closed   subspace  of
$L^m_{t^-_k,\delta}(f^*N_V)$ consisting of all sections  that vanish
to order $s_{k,i}$ at $p_{k,i}$ and order $t^+_{k,i}$ at
$y^+_{k,i}$. By standard elliptic theory for weighted norms
(cf.\ \cite{L}) the operator
     $D^N$ defines a bounded operator
\bear
D^N_k: L^m_{k,\delta}(f^*N_V) \to L^{m-1}_{k,\delta+1}(T^*C\otimes f^*N_V)
\label{7.Lspace}
\eear
which, for generic $0<\delta<1$, is Fredholm with
$$
           \mbox{index}_\R\  D_{k}^N\ =\
           2c_1(N_V)A_k+\chi(B_k)+2(\deg t_k^- -\deg s_k -\deg t_k^+)=\chi(B_k)
\label{indexDNlayerk}
$$
where $A_k=[f(B_k)]$ in $H_2(X)$.  We used the fact that
$c_1(N_V)A_k= \deg s_k +\deg t_k^+-\deg t_k^-$ (since the Euler class
of a line bundle can be computed from the zeros and poles of a
section). Lemma \ref{DNcxlemma} implies that the kernel of this operator is
$J$-invariant, and so we can form the complex projective space $\P(\ker
D_k^N)$.

\numbereddemo{Definition} \label{VstableSMap} A {\it $V$\/{\rm -}\/stable map} is a stable map
$(f,\phi)\in\ov\M_{g,n+\ell(s)}(X,A)$  together with
\begin{itemize}
\item[(a)] a marked layer structure on its domain $C$ with
$f|_{B_0}$ being $V$-regular, and

\item[(b)]   for each $k\ge 1$ an  element  $[\xi_k]$ of
$\P(\ker D_k^N)$  defined on the layer $B_k$ by a section $\xi_k$
that is nontrivial on every irreducible component of $B_k$.
\end{itemize}
\enddemo

Let $\ov\M^V_{g,n,s}(X,A)$ denote the set of all $V$-stable maps.
This contains the set
$\M^V_{g,n,s}(X,A)$  of $V$-regular maps as the open subset --- the
$V$-stable maps whose entire domain lies in  layer 0.  Forgetting the
data  $[\xi_k]$ defines a map
\bear
\ov\M^V_{g,n,s}(X,A)\ {\displaystyle\mathop\rightarrow_{\beta}}\
\ov\M_{g,n+\ell(s)}(X,A).
\label{7.spaceofVstable}
\eear

Each $V$-stable map $(f, \phi,[\xi_1],\dots, [\xi_r])$ determines an
element of the space ${\cal H}^V_X$ of Definition \ref{5.defcalH} as
follows.  For a very small $\ep$, we can push the components in $V$
off $V$ by composing $f$ with $\mbox{exp}(\ep^k \xi_k)$ and,
for each $k$, smoothing the domain at the nodes
$B_k\cap\l(\ma\cup_{l>k}B_{l}\r)$ and smoothly joining the images
where the zeros of $\ep^k\xi_k$ on $B_k$ approximate the poles of
$\ep^{k+1}\xi_{k+1}$.  The resulting map
$$
f_\xi\,=\,f|_{B_0} \,\#\,  \exp(\ep\xi_1)\,\#\,\cdots \,\#\, \exp(\ep^r\xi_r)
$$
is $V$-regular, and so represents a homology class
$h(f,\phi,[\xi_k])=h(f_\xi)\in{\cal H}^V_X$
under (\ref{5.defh}). That class depends only on $[\xi]$:  for
different choices of the $\ep_k$ and of
representatives of the $[\xi_k]$, the $\ep^k\xi_k$  are homotopic
through nonzero elements of the
kernel with the same zeros and poles and hence represent the same element of
${\cal H}^V_X$. Thus there is
a well-defined map
\bear
\ov\M^V_{g,n,s}(X,A)\; \ma\longrightarrow^h\ {\cal H}^V_{X,A,s}.
\label{7.mapdh}
\eear

\proclaim{Proposition}  \label{7.Vspecial2}
There exists a topology on   $\ov\M^V_{g,n,s}(X,A)$
which makes it compact and for which
the maps $\beta$ of {\rm (\ref{7.spaceofVstable})} and $h$ of
{\rm (\ref{7.mapdh})}  are continuous and differentiable on each stratum.
\endproclaim 

\pf  There are three steps to the proof.  The first looks at
sequences of $V$-regular maps (which are
$V$-stable maps with trivial layer structure) and the  second
analyzes a general sequence  of
$V$-stable maps.  The third step uses that analysis to define the
topology on $\ov\M^V_{g,n,s}(X,A)$.

      Let $f_m:C_m\to X$ be a  sequence of maps in $\M^V_{g,n,s}(X,A)$.  By
the bubble tree convergence Theorem
\ref{BTCthm}, a subsequence, still called $f_m$,  converges to a
stable map $f:C\to X$. By successive renormalizations we will give
the limit map $f$ the structure of
a    $V$-stable map $(f,[\xi])$.

Since the last $\ell(s)$ marked points converge, the multiplicity
vector $s$ of $f_m$ carries over to the limit, defining the vector $s$
of Definition \ref{7.markedlayerstructure}b.  The rest of the layered
structure
is defined inductively.  We assign $\lambda_j=0$ to each component
$C_j$ unless $f(C_j) \subset V$, so that the layer $B_0$ consists of all
components that are not mapped into $V$.  Let $C(1)$ be the union of
those components of $C$ not in layer 0.  Assign to each double point
$y$ of $B_0\cap C(1)$ a multiplicity $t_y$ equal to the order of
contact of $f|_{B_0}$ with $V$ at $y$.

Now apply the  argument of Proposition
\ref{TaubesRenomlemma2}.  That  produces  renormalized maps
$g_{m,1}=\exp \xi_{m,1}$ which converge to
a  nontrivial element of $\xi_1$ in $\ker D^N$ on  $C(1)$.   We
assign $\lambda_j=1$ to
each component $C_j\subset C(1)$ on which $\xi_1$ is nonzero and
denote the union of the remaining
components by $C(2)$. Then $\xi_1$ is defined and nonzero on every
component of $B_1$.  Moreover,

\begin{itemize}
\item[(a)]  $\xi_1$ vanishes at the double points $y$ where $B_1$ meets
$C(2)$.  We assign such $y$  a multiplicity $t_y$ equal to the order
of vanishing of $\xi_1$ at that point.

\item[(b)] As in the proof of Proposition \ref{TaubesRenomlemma2}, $\xi_1$
has a pole of order $t_x$ at each $x\in B_0\cap B_1$ and vanishes to
order $s_{1i}$ at the points $p_{1i}$ in $B(1)$.

\item[(c)] $\xi_{m,1}\to \xi_{1}$ and hence $f_m$ and
$f|_{B_0} \,\#\,\exp(\ep\xi_1)\,\#\, g_{m,1}|_{C(2)}$
define the same element of ${\cal H}^V_X$ for large $m$.
\end{itemize}
This defines $[\xi_1]$ and multiplicity vectors $s$ and
$t$  on $B_1$.

Next, repeat the renormalization on $C(2)$ and continue.  This
inductively defines a layer structure on $C$,  multiplicities $t$ for each
point in $B_k\cap C(k+1)$,
and determines  a nontrivial element $\xi_k$ of  the kernel of
$D^N_k$ on each layer $B_k$.  This process terminates because
each $C(k)$ has fewer components
than $C(k-1)$ while parts (a) and (c) of Lemma \ref{minenergylemma}
give a uniform bound on the
number of components. The end result is a nontrivial $\xi_k$
on every component of $B_k$, $k\geq 1$. From (c) above we see  that
$\lim [f_m]=h(f,[\xi_k])$ in ${\cal
H}^V_{X,A,s}$.

\vglue3pt

      We next consider general sequences in $\ov\M^V_{g,n,s}(X,A)$.
Given a sequence $F_m=(f_m, [\xi_{m,k}])$ of $V$-stable maps, we first
form the maps $f_{m,0}$ obtained by restricting $f_m$ to its bottom
layer $B_m(0)$.  The $f_{m,0}$ may represent different homology
classes, but there is a uniform bound on their energy, so that the iterated
renormalization argument of Lemma \ref{7.Vspecial2} produces a
subsequence converging to a $V$-stable map $F_0$.  Similarly, the
restrictions of $f_m$ to $B(1)$ converge to a stable map $f_1$ into
$V$ and the renormalized maps $\exp(\xi_{m,1}):B_m(1)\to \P_V$ have a
subsequence converging to a limit which, on its bottom layer, has the form
$\exp (\xi_{0,1}$) with $\xi_{0,1}\in\mbox{ker } D^N$.  Then
$F_1=(f_{m,1}, [\xi_{0,1}])$ is a $V$-stable map whose image lies in
$V$ and whose bottom layer fits with the top layer of $F_0$ to form a
$V$-stable map $F_0\cup F_1$. This process continues, and terminates
because   each layer  carries energy at least $\alpha_V$.

Finally, observe that this renormalization process can be read
differently: it actually {\it defines} a notion of a convergence
sequence of $V$-stable maps.  Convergence in that sense defines a
topology on $\ov\M^V_{g,n,s}(X,A)$, which we adopt as {\it the}
topology on the space of $V$-stable maps.  Reinterpreted, the above
analysis shows that $\ov\M^V_{g,n,s}(X,A)$ is compact and $h$ is
continuous with that topology.  \enddemo

      The next theorem is the key result needed to define the relative
invariants; it implies and supersedes
Proposition \ref{compactificationThm}.

\proclaim{Theorem}  \label{7.Vstablemapsiscompact}
 The space of $V$\/{\rm -}\/stable maps is compact and there is
a continuous map
\bear
\ep_V:\ov{\M}_{g,n,s}^V(X,A) \lraop{{\rm {st}  \times {\rm{ev}}\times h}}
\
\ \ov{\M}_{g,n+\ell(s)}\times X^n
\ti {\cal H}^V_{X,A,s}
\label{7.MainTheorem}
\eear
obtained from {\rm (\ref{2.MainEvalmap})} and {\rm (\ref{7.mapdh}).}
Furthermore{\rm ,} ${\M}_{g,n,s}^V(X,A)$ is oriented and the complement of
      ${\M}_{g,n,s}^V(X,A)$ in the irreducible part of
$\ov{\M}_{g,n,s}^V(X,A)$ has codimension  at least two.
\endproclaim

\pf   To define the orientation, note that at each $f\in \M^V_s(X)$, the
tangent space is the kernel of the linearized operator $D_f$.  For
generic $(J,\nu)$ the cokernel vanishes, so the tangent space is
identified with the formal vector space
$\mbox{ker}\,D_f-\mbox{coker}\,D_f$.  This is oriented by the mod 2
spectral flow of a path in the space of Fredholm operators that
connects $D_f$ to any operator that commutes with $J$, where the
kernel and the cokernel are complex vector spaces and hence are
canonically oriented.  This orients $\M^V_s(X)$, and the orientation
extends to the compactification provided the frontier strata have
codimension at least 2.

       We established compactness above and the dimension statements are
verified in the next two lemmas.
      \enddemo

For the dimension counts we return to the notation of
(\ref{7.Lspace}).  The strata of the space of $V$-stable maps are
labeled by curves $C=\bigcup B_k$ with a marked layer structure but no
specified complex structure; we denote these strata by
$\ov\M^V_{g,n,s}(C)$.  For each layer $B_k$, $k\ge 1$, let
\bear
\ov\M^V_{B_k,n_k,s_k\cup t^+_k,t^-_k}(X,A_k)
\label{7.strata}
\eear
denote the set of $V$-stable maps $(f,[\xi])$ where
$f\in\M_{B_k,n_k+\ell(s_k\cup t^+_k\cup t^-_k)}(V,A_k)$ is a map from
$B_k$ to $V$ for which
moreover the operator (\ref{7.Lspace}) has a nontrivial kernel on
each irreducible component of $B_k$.

In this context consider the  Hilbert bundle over the universal moduli space
$$
L^M_{k}(V) \to \ov\UM_{B_k,n_k+\ell(s_k\cup t^+_k\cup t^-_k)}(V,A_k)
$$
    whose fiber at $f:B_k\to V$ is the space  $L^{m}_{k,\delta}(f^*N_V)$ of
(\ref{7.Lspace}).  It is straightforward to adapt the proof of Lemma
\ref{transverseLemma1} to show that $D^N_k$ defines a
section of
\bear
\begin{array}{c}
{\rm Fred}\l(L^m_{s,\delta}(V), L^{m-1}_{s,\delta+1}(V)\r)
    \\[5pt]
\big\downarrow\\[5pt]
\ov\UM_{B_k,n_k+\ell(s_k\cup t^+_k\cup t^-_k)}(V,A_k)
\end{array}
\label{7.diagramFred}
\eear
which is transverse to the subspaces ${\rm Fred}_r$ of operators
with kernel of dimension $r\geq 1$.

\proclaim{Lemma} The irreducible part of the  space {\rm (\ref{7.strata})}
is an orbifold of {\rm ``}\/correct\/{\rm ''} dimension{\rm ,} which is at most
\bear \qquad\quad
d_k&\hskip-7pt=\hskip-7pt&2\left[-K_X[A_k]+\frac14 (\dim X-6)\chi(B_k)+n_k+\ell(t_k^-) +
\ell(s_k)+\ell(t_k^+)\r.
\label{6.dimM'}\\
&\hskip-7pt\hskip-7pt&\l.\hskip.6in +\deg t_k^--\deg t_k^+-\deg s_k \r]-2.
\nonumber
\eear
\endproclaim

\pf Consider an irreducible  $V$-stable map  $(f,[\xi])$ in
the space  (\ref{7.strata}). By Theorem \ref{BTStructurethm}b,
generically $\M_{B_k,n_k+\ell(s_k\cup t^+_k\cup t^-_k)}(V,A_k)^*$ is
an orbifold of dimension
\begin{equation}
d_k'=2\l[-K_V[A_k]+\frac14 (\dim V-6)\chi(B_k) +n_k+\ell(t_k^+)+\ell(t_k^-)+
\ell(s_k)\right].\enspace
\label{6.dimM(V)}
\end{equation}
Comparing (\ref{6.dimM'}), (\ref{6.dimM(V)}) and the displayed formula following
(7.1) we see that
$d_k-d_k'=2(\iota-1)$  where  $2\iota=\mbox{index }D^N_k$. The lemma
follows immediately if $\iota >0$.

When $\iota \le 0$  we can
use the element $\xi_k$  and the tranversality of $D^N_k$ in
(\ref{7.diagramFred}) to conclude that
the set of $f\in\M_{B_k}(V,A_k)^*$ with $\mbox{dim }\ker D^N_k= 2r$ form a
suborbifold of codimension $2r(-\iota+r)$.
The lemma then  follows because
$r(r-\iota)\ge 1-\iota$ for all $r\ge 1$.
(Notice that this argument requires that $\xi_k$ be nontrivial on every
component of $B_k$ as in Definition \ref{VstableSMap}). \enddemo

\proclaim{Lemma} Each irreducible stratum $\ov\M^V_{g,n,s}(C)$ is an
orbifold whose dimension is $2(r+\sum \ell_k)$  less than that in  {\rm (\ref{5.1}),}
 where $r$ is the total number of nontrivial layers and
$\ell_k$ is the number of double points of $C$ in layer $k\ge 0$. 
\endproclaim

\pf  The stratum $\ov\M^V_{g,n,s}(C)$ is the product of the spaces
(\ref{7.strata}), one for each  layer, constrained by the matching
conditions $f(x)=f(y)$ at each of the $\ell(t)$ double points where $B_k$
meets the other layers.  A standard tranversality
argument \cite{RT1} shows that the irreducible part of this space is
an orbifold of the expected dimension.
Thus
$$
\dim  \ov\M^V_{g,n,s}(C)\ \le \ \sum_{k=0}^r d_k\ -\ \ell(t)\  \dim \, V.
$$
Now substitute in  (\ref{6.dimM'}) for  $d_k$, $k\ge 1$ and sum, noting that
   (i) formula (\ref{6.dimM'}) is additive in  $A_k$ and  $s_k$,
   (ii)  each double point $x$ between different layers contributes its
multiplicity $t_x$ to both $t^+$ and
$t^-$ so that  $\sum \deg t_k^+ = \sum \deg t_k^-$, and
   (iii) the Euler characteristics
add according to  the formula
\best
\sum \chi(B_k) =\chi(C) +2\ell(t)
\eest
where $\ell(t)=\sum \ell(t_k^+)=\sum \ell(t_k^-)$.  The result
follows. \enddemo

\numbereddemo{{R}emark} There is a different but equivalent viewpoint on what a
$V$-stable map is.  In Proposition
\ref{TaubesRenomlemma2} and Lemma \ref{7.Vspecial2} we inductively
produced  a limit map $f$, layers $B_k$, and limiting renormalized
maps $g_k:B_k\to \P_V$ for $k\geq 1$.  On each $B_k$ we  wrote $g_k$ as
$\exp_{f}(\xi_k)$ using the exponential map from
$V_0\subset P_V$; the information $\{[\xi_k]\}$ then
defined a  $V$-stable map as in
Definition \ref{VstableSMap}.  Alternatively, we could have recorded
the $g_k$ themselves modulo the $\cx^*$
action on $\P_V$.  From that perspective the limiting $V$-stable map is an
equivalence  class of continuous maps
$f\cup g_1\cup
\cdots \cup g_r$ from
$C=\ma\bigcup_{k=0}^r B_k$  into the singular space
\best
X\ma\cup_{V=V_\infty}\P_V \ma \cup_{V_0=V_\infty} \cdots \ma
\cup_{V_0=V_\infty} \P_V
\eest
with $f$ mapping $B_0$ to $X$ and each $g_k$, $k\ge 1$, mapping $B_k$
to the $k^{\rm th}$ copy
of $\P_V$ with all maps $V$-regular along each intermediate copy of
$V$ and   with two \pagebreak such maps $f\cup g_1\cup
\cdots \cup g_r$ equivalent  if they lie in the
same orbit of the $(\cx^*)^r=\cx^*\times\cdots\times \cx^*$ action. The
renormalization procedure of  Lemma
\ref{7.Vspecial2} gives a uniform bound on the number of  copies of $\P_V$  for
each homology class~$A$.

The correspondence between these two viewpoints is clear.  We found
that the  analytic technicalities
were  easiest using the description of Definition \ref{VstableSMap},
but in general it is useful to keep
both descriptions in mind.
\enddemo

\section{Relative invariants}

The relative Gromov-Witten invariant is the homology class obtained by
pushing forward the compactified space of $V$-regular maps by
the map (\ref{7.MainTheorem}):
$$
\ep_V:\ov\M_{g,n,s}^V(X,A) \ \longrightarrow\    \ov\M_{g,n+\ell(s)}\times X^n
\ti {\cal H}_X^V.
$$
In this section we show that this yields a well-defined homology class. We
then recast the relative
invariants as Laurent series and explain their geometric interpretation.

\proclaim{Theorem} Assume $X$ and $V$ are semipositive or that the moduli
space
$\ov\M_{g,n,s}^V(X,A)$ is generically irreducible. Then for generic
$(J,\nu)\in {\cal J}^V${\rm ,}  the image of $\ov\M_{g,n,s}^V(X,A)$ under
$\ep_V$ defines an element
\bear
{\rm GW}^V_{X,A,g,n,s} \in H_*(\ov\M_{g,n+\ell(s)}\ti X^n\ti {\cal H}_X^V ;\Q)
\label{7.defclass}
\eear
of dimension
\bear
-2K_X[A]+(\dim X-6)(1-g)+2(n+\ell(s) - \deg s).
\label{index.0}
\eear
This homology class is independent of  the generic $(J,\nu)\in {\cal J}^V$.
\label{welldefn}
\endproclaim

\pf  The spaces ${\cal H}_X^V$ and $\ov\M_{g,n}$ are orbifolds, and
therefore so is $Y=\ov{\M}_{g,n}\ti X^n \ti {\cal H}_X^V$.
Fix a generic $(J,\nu)$ and consider the image of the smooth map
$\ep_V:\M^V_{g,n,s}\to Y$.  Its frontier
$$
{\rm Fr}(\ep_V)\ =\ \left\{ y\in Y\ |\ y=\lim\;\ep_V(f_k) {\ \rm and\ }
\{f_k\}\  \mbox{\rm has  no
convergent subsequence}\ \right\}
$$
is exactly the image
$$
\left[\ep_V\left(\ov{\M}_{g,n,s}^V(X,A)\setminus
{\M}_{g,n,s}^V(X,A) \right)\right].
$$
Then Theorem \ref{7.Vstablemapsiscompact} (applied to the reduced
moduli space when $X$, $V$ are semipositive) implies that the frontier
${\rm Fr}(\ep_V)$ lies in a set of dimension two less that the dimension
(\ref{index.0}) of the image.  Proposition 4.2 of \cite{KM} then
implies that the image
$\ep_V\left(\ov{\M}^V_{g,n,s}\right)$ carries a rational homology class
of that dimension.

The last statement of the theorem follows by a cobordism argument.
By  Theorem
\ref{JConnectedtheorem} of the appendix ${\cal J}^V$ is path-connected; so
any two generic pairs
$(J_1,\nu_1)$ and $(J_2,\nu_2)$ can be joined by a path $\gamma(t)$.    The
Sard-Smale theorem, applied to
      the space of such paths, shows that over the generic such $\gamma(t)$, the
universal moduli space $\UM^V$ over $\gamma$ is an orbifold.
Again the frontier of the image
$$
\ep_V\left(\pi^{-1}(\gamma)\right)
$$
lies in a set of dimension two less that the dimension  of this image.
Proposition~4.4 of \cite{KM} then
implies that the homology classes (\ref{7.defclass}) defined by
$(J_1,\nu_1)$ and $(J_2,\nu_2)$ are the
same.
\enddemo

\numbereddemo{Definition} Let $(X,\w)$ be a closed symplectic manifold with a
codimension two symplectic
submanifold $V$.   For each $g,n$, the  {\it relative {\rm GW} invariant} of
$(X,V,\w)$ is the
homology class (\ref{7.defclass}).
\enddemo

It is again convenient to assemble these invariants into a  Laurent
series.
For that, we simply repeat the
discussion leading from (\ref{modulimap2}) to (\ref{defGW}). Thus the full
relative GW invariant is  the
map
\bear
{\rm GW}_X^V: H^*(\ov\M)\otimes  {\Bbb T}^*(X) \ \longrightarrow\
      H_*({\cal H}_X^V;\Q[\la])
\label{relativeGW}
\eear
where ${\Bbb T}^*(X)$ is the total tensor algebra ${\Bbb T}(H^*(X))$ on the
rational cohomology.
\smallskip

As in (\ref{5.defh2}), ${\cal H}^V_X$  is a union of components labeled by
$A$ and $s$, so that
$$
H_*\left({\cal H}_X^V;\Q[\la]\right)\ =\ \ma \bigoplus_{A,s\atop
\mbox{\tiny deg } s =A\cdot V}\
H_*\left({\cal H}^V_{X,A,s}\right)\otimes \Q[\la].
$$
Thus there is an expansion
\bear
{\rm GW}_X^V =\ \sum_{g,n} {1\over n!}\ \ma \sum_{A, s\atop \mbox{\tiny deg }
s=A\cdot V} {1\over \ell(s)!}\; {\rm GW}^V_{X,A, g,n,s}\ t_A\ \la^{2g-2} 
\label{defnRelInvt2}
\eear
where the coefficients on the right lie in
$H_*\left({\cal H}^V_{X,A,s}\right)$.

Formula (\ref{defnGT}) extends this to a  relative Gromov-Taubes invariant
\begin{equation}
{\rm GT}_X^V= \mbox{exp}\,({\rm GW}^V_X):  H^*(\wt\M)\otimes  {\Bbb T}^*(X)
\to    H_*({\cal H}^V_X;\Q[\la]).
\label{defn5}
\end{equation}

It is clear that these invariants are natural:  if $\phi$   is a
diffeomorphism of $X$ then
$V'=\phi^{-1}(V)$ is a symplectic submanifold of $(X',\phi^*\w)$ and
$$
{\rm GW}^{V'}_{X'}\ =\ {\rm GW}^V_X.
$$

It is also  clear that these invariants extend the GW invariants of Section~2.  In fact, the entire
construction carries through when  $V$ is  the empty set.  In that case
${\cal H}_X^V$ is just $H_2(X)$
and the relative invariant takes values
in $NH_2(X)$.  The relative and absolute invariants are then equal:
$$
{\rm GW}_X^{\emptyset}={\rm GW}_X.
$$

More importantly, the relative GW invariants are unchanged under  symplectic
isotopies, i.e. they are constant as we move along
1-parameter families $(X, V_t,\w_t)$  consisting of   a symplectic form
$\w_t$ on $X$ and a codimension-two  submanifold $V_t$  which is symplectic
for $\w_t$.
More generally,  we say
$(X,V,\w)$  is {\it deformation equivalent} to  $(X',V',\w')$ if there is a
diffeomorphism $\phi:X'\to X$ such that $(X', \phi^{-1}(V),\phi^*\w)$ is
isotopic to $(X',V',\w')$

\proclaim{Proposition}
\label{7.invarianceofGW}
The relative invariant ${\rm GW}^V_X$ depends only on the symplectic deformation
class of $(X, V,\w)$.
\endproclaim 

\pf  By naturality we need only verify invariance under symplectic
isotopies. But that follows by
essentially the same cobordism argument used in the last part of the proof
of Theorem \ref{welldefn}  (cf.\
\cite[ Lemma 4.9]{RT2} and \cite[Prop.\ 2.3]{LT}).
\enddemo

The geometric meaning of the relative invariant is obtained by evaluating
the homology classes
(\ref{7.defclass}) on dual cohomology classes, thereby re-expressing the
invariant as a collection of
numbers.  To do that, choose $\kappa\in H^*(\ov\M_{g,n})$, a vector
$\alpha=(\alpha_1,\dots, \alpha_n)$ of classes in $H^*(X)$, and a
$\gamma\in H^*({\cal H}_X^V)$, such that
\bear
      \deg \kappa+\deg \al -2\ell(\al)+\deg \gamma
\label{7.dimrestriction}
\eear
      is the dimension (\ref{index.0}) of the homology class ${\rm GW}^V_{X, 
A,g,n,s}$.
Then the evaluation pairing
gives numbers
\bear
{\rm GW}^V_{X,A,g,n,s}(\kappa,\alpha,\gamma)\ =\ \left<\left[{\rm GW}^V_{X,A,g,n,s}\right],
\kappa\cup\alpha\cup\gamma
\right>.
\label{7.pairing}
\eear
It is these numbers that have a specific geometric interpretation.

\proclaim{Proposition} \label{prop7.2}
Fix generic geometric representatives $K\subset \overline{\M}_{g,n}${\rm ,}
$A_i\subset X${\rm ,} and
$\Gamma\subset {\cal H}_X^V$  of the  Poincar{\rm \'{\it e}}
duals of $\kappa${\rm ,} $\alpha$ and $\gamma$. Then the evaluation
{\rm (\ref{7.pairing})} counts   the oriented number of
genus $g$ $(J,\nu)$\/{\rm -}\/holomorphic maps\break $f:C \to X$ with $C\in K${\rm ,}
$h(f)\in \Gamma${\rm ,} and $f(x_i)\in A_i$ for each marked point $x_i$.
\endproclaim

Note that the condition $h(f)\in \Gamma$, defined by the map
(\ref{5.defh}), constrains both the homology
class $A$ of the map and the boundary values of the curve.  In the special
case when there are no rim tori,  these homology and the boundary value
constraints can be fully separated as  in the beginning of Section 9.

Also note that the invariant counts maps from a domain with $n+\ell(s)$
marked points, the last $\ell(s)$ of which are mapped into
$V$.  Two such maps with their last $\ell(s)$ marked points renumbered
are considered different.  This might seem to introduce an unnecessary
redundancy, but the marking on the last set of points is  needed to prove
that  two curves whose intersection with $V$ is the same can be `glued
together' (see \cite[Th.\ 5.6]{IP4}).

It is useful to have a more general version of Proposition \ref{prop7.2}.
If we drop the dimension
restriction (\ref{7.dimrestriction}) then the set of $(J,\nu)$-holomorphic
maps
\bear
\M_{g,n,s}^V(X,A;\kappa,\al,\gamma)
\label{7.cutdownMV}
\eear
satisfying the conditions of Proposition
\ref{prop7.2} will no longer be finite; its expected dimension is
\bear
\label{7.dimMkag}
\dim \M_{g,n,s}^V(X,A)-\deg \kappa-\deg \al +2\ell(\al)-\deg\gamma
\eear
with the first term given by (\ref{index.0}). Of course (\ref{7.cutdownMV})
depends on $(J,\nu)$; it is
a fiber  of a map
$$
\UM_{g,n,s}^V(X,A;\kappa,\al,\gamma) \to {\cal J}^V
$$
from the subset of the universal space $\UM_{g,n,s}^V(X,A)$  consisting of
those maps satisfying the
conditions of Proposition \ref{prop7.2}.  Proposition \ref{prop7.2} is then
the 0-dimensional case of the following fact.

\proclaim{Lemma}  \label{Lemma6.2}
For generic representatives $K, A_i$ and $\Gamma${\rm ,}  
$$\UM_{g,n,s}^V(X,A;\kappa,\al,\gamma)^*$$ is an orbifold.   Hence for generic
$(J,\nu,K,A_i,\Gamma)${\rm ,} the irreducible part of\break
$\M_{g,n,s}^V(X,A;\kappa,\al,\gamma)$  is an orbifold of dimension as in
{\rm (\ref{7.dimMkag}).}
\endproclaim

\pf Consider the evaluation map ${\mbox{st}}\ti {\mbox{ev}}:
{\UM}^V_{g,n,s}(X,A)\ra \M_{g,n}\ti X^n\ti V^{\ell(s)}$ by
\begin{eqnarray*}
&&(f,j,x_1, \dots, x_n, (p_1,s_1)\dots,(p_\ell,s_\ell))\\
&&\qquad\quad\mapsto
((j,x_1,\dots,x_n),f(x_1), \dots,f(x_n),f(p_1),\dots,f(p_\ell)).
\end{eqnarray*}
We can achieve tranversality by moving the geometric representatives
$K$, $A_i, C_l$ of $\al,\gamma$, exactly as in Section 4 of
\cite{RT1}, only now keeping the representatives $\Gamma$ in $V$.  The
dimension formula (\ref{7.dimMkag}) follows because $\kappa$ cuts down
$\deg \kappa$ dimensions, and the $\alpha$ and $\gamma$ constraints
cut down by $\deg(\al)$ and $\deg(\gamma)$ dimensions,
respectively. \enddemo

\numbereddemo{{R}emark} So far, the relative invariant is defined by cutting down
with geometric representatives of classes in $H_*(X)$.  For the
application to the gluing theorem in \cite{IP4} it is useful to have a
version of these invariants that allows constraints in
$H_*(\hat{X},S)$ (which are Poincar\'e dual to classes in
$H^*(X\setminus V)$).

      Let $\hat{X}$ be the manifold (\ref{hatX}) obtained from $X\setminus
V$ by attaching as boundary a copy of the unit circle bundle $\pi:S\to
V$ of the normal bundle of $V$ in $X$.  Suppose that $Z$ is a
symplectic sum obtained by gluing $\hat{X}$ to a similar manifold
$\hat{Y}$ along $S$.  We can then consider stable maps in $Z$
constrained by classes $B$ in $H_k(Z)$, i.e. the set of stable maps
$f$ with the image $f(x)$ of a marked point lying on a geometric
representative $\phi$ of $B$.  If we restrict ourselves to the $\hat{X}$ side,
$\phi$ defines a class $[\phi]\in H_*(\hat{X},S)$ and a subspace
$\ov\M^V(\phi)$ of the space of $V$-stable maps.

We can repeat the above arguments to obtain GW invariants constrained
by $\phi$.  However, in general $\ov\M^V(\phi)$ is an orbifold with
boundary and consequently the invariant depends on the choice of the
representative $\phi$.  In \cite[\S13]{IP4} we will show that
this dependence is rather mild and in practice can be handled as
follows.

In the exact sequence
$$
H_k(X\setminus V)=H_k(\hat{X})\longrightarrow H_k(\hat{X},S)
\stackrel {\bd}\longrightarrow  H_{k-1}(S)
\stackrel {\iota}\longrightarrow H_{k-1}(\hat{X})
$$
choose a splitting at $\bd$ over $K=\mbox{ker }\iota$ and choose a
geometric representative $\phi_B$ of each
$B$ in the image of this splitting. For each class in $H_k(\hat X, S)$ we
can
find a geometric representative $\phi$  with the same
boundary as one of the chosen $\phi_B$.  Then $a=[\phi\#(-\phi_B)]$ lies in
$H_k(X\setminus V)$ and
\bear
GW(\phi)=GW(\phi_B)+GW(a).
\label{8.additive}
\eear
This defines ${\rm GW}$ on a complete set of  representatives of $H_k(X;V)$
which combines with Poincar\'e duality  to give a  map
\bear
{\rm GW}:H^*(\ov\M)\ti {\Bbb T}^* (X\setminus V)\longrightarrow H_*({\cal H}^V_X).
\eear
This map is neither canonical nor linear in the constraints, but is
additive as in (\ref{8.additive}) and is
determined by the invariants (\ref{relativeGW}) and those for the chosen
constraints $\phi_B$.  It provides a set of constraints that can be used in
the gluing theorem to constrain by any class in $H^*(Z)$.
\enddemo

\section{Examples}

The relative GW invariants are designed to be used in `cutting and
pasting' arguments of symplectic topology.  However, in several
interesting cases they are identical to invariants from enumerative
algebraic geometry.  Here we give three such examples; in each case
$X$, $V$ are semipositive.  Actual computations of the relative
invariants in these cases are done in \cite{IP4}.

As noted in Remark 5.3, the description of the relative GW invariants is
simplified considerably when there are no rim tori.  This occurs
whenever $H_1(V)=0$ and more generally when every rim torus represents
zero in\break $H_2(X\setminus V)$.  In these cases there is no covering
(\ref{5.HSVproj}), and ${\cal H}_X^V$ is the subset of $H_2(X) \ti
{\cal S}V$ consisting of pairs $(A,s)$ with $\deg s=A\cdot V$.  The
homology of ${\cal H}_X^V$ is the corresponding subalgebra of
$NH_2(X)\otimes {\Bbb {CT}}_*(V)$ where ${\Bbb {CT}}_*(V)$ is the ``contact
tensor algebra'' of $V$:
\best
{\Bbb {CT}}_*(V)\,=\, {\Bbb {T}}({\Bbb N}\ti H_*(V)).
\eest
The relative invariants are then maps
\best H^*(\ov\M)\otimes {\Bbb T}^*(X) \ra {\Bbb {CT}}_*(V)\otimes
NH_2(X)[\la]
\eest
and have Laurent expansions like
(\ref{defnRelInvt2}) with coefficients in ${\Bbb {CT}}_*(V)$, and
$NH_2(X)$ is the Novikov ring (cf.\ Section 2).  Fix dual bases
$\gamma_i$ of $H_*(V;\Q)$ and $\gamma^i$ of $H^*(V;\Q)$.  Then bases
of the contact algebra and its dual are given by elements of the form
\begin{equation}
C_{s,\gamma}=C_{s_1, \gamma_1}\otimes\dots \otimes C_{s_\ell,
\gamma_\ell}
\quad \mbox{ and }\quad
C_{s,\gamma}^*=C_{s_\ell, \gamma^\ell}\otimes\dots \otimes C_{s_1,
\gamma^1}\qquad\quad
\label{last.basis}
\end{equation}
respectively, where  $s_i\ge 1$ are integers.
With $\kappa$ and $\alpha$ as above, we can expand
\bear
{\rm GW}^V(\kappa,\alpha)= \sum_{A,s,\gamma}{1\over \ell(s)!} \;
{\rm GW}_{X,A,g}^V(\kappa,\alpha; C_{s,\gamma}^*)\ \
C_{s,\gamma}\ t_A\ \la^{2g-2}
\label{2.last}
\eear
where the coefficients count the oriented number of genus $g$
$(J,\nu)$-holomorphic, $V$-stable maps $f:C \to X$ with $C\in K$,
$f(x_i)\in A_i$; these have a contact of order $s_j$ with $V$ along
fixed representatives $\Gamma_j$ of the Poincar\'e duals of the
$\gamma^j$, where $K$ and $A_i$ are Poincar\'e duals of the $\kappa$ and
$\al_i$.

With this background, we will describe two simple examples of relative
invariants and their classical algebraic-geometry counterparts.  While
these are amongst the very simplest examples of relative invariants,
each has a long history and has proved to be frustratingly difficult
to compute by algebraic-geometric methods.  However, in both examples
recent progress has been made on calculating the invariants by using,
in part, symplectic cut-and-paste arguments.

\numbereddemo{Example} \label{ex9.1}
 The Hurwitz numbers are examples of GW invariants of
$\P^1 $ relative  to several points in $\P^1$.
\enddemo

The classical Hurwitz number $N_{g,d}$ counts the number of
nonsingular, genus $g$ curves expressible as $d$-sheeted covers of
$\P^1$ with a {\it fixed} branch divisor in general position. They
were first computed in \cite{Hu} by combinatorial techniques. More
generally, if $\al$ is an (unordered) partition of $d$ then the
Hurwitz number $N_{g,d}(\al)$ counts the number of smooth degree $d$
maps from a genus~$g$ Riemann surface to $\P^1$ with the ramification
above a fixed point $p_0$ as specified by the partition $\al$, and
simple branching at exactly $r(g,\al)=d+\ell(\al)+2g-2$ other fixed
points in general position.

On the other hand, for any distinct fixed points $p_0,\dots, p_{r}$ in
$X= \P^1$, the set $V=\{p_0,\dots, p_{r}\}$ is a symplectic
submanifold of $X$ with no rim tori.  The homology class $A$ of the
map is given by the degree $d$; we take $\kappa=1$ because we are
imposing no constraints on the complex structure of the curves.  Thus
the relative GW invariant ${\rm GW}_g^{V}(\P^1,d)$ has the form
(\ref{2.last}) with values in ${\cal S}V$.  But ${\cal S}\{p_0\}$ is
the disjoint union of copies of $p_0$, one copy for each vector $s$
with $\mbox{deg}\ s=p_0\cdot A=d$. Furthermore, the relative invariant
is unchanged when the marked points that are mapped into $p_0$ are
permuted.  We can then associate the generator of $H_0(V_s)$ coming
from the point $p_0$ with the monomial $z^s=z_{s_1}\cdots
z_{s_\ell}$. Thus we identify
$$
H_*({\cal S}\{p_0\})\ =\ \bigoplus_s\ \Z\ = \ \Z[z_1,z_2,\dots]
$$
where the last term is the polynomial ring on variables
$z_1,z_2,\dots\, $. Then the  Hurwitz number
\best
N_{g,d}(\al)\ =\  {\rm GW}_{\P^1,d,g}^{V} (z^{s};b^r)
\eest
where $r=d+2g-2+\ell(\al)$ and $s$ is one -- any one -- of the ordered
partitions obtained by ordering $\al$.  Geometrically, the variable
$z_i$ models a contact of order $i$ at $p_0$, $r$ is the number of
leftover simple branch points, and $b^r$ denotes the condition that
$r$ simple branch points are mapped to $r$ distinct fixed points
$\{p_1,\dots,p_r\}$.

With this notation the Laurent series (\ref{defnRelInvt2}) of the relative
invariant is:
\begin{eqnarray}
{\rm GW}^p_{\P^1}& =&\sum_s \ {1\over \ell(s)!} {\rm GW}_{\P^1,d,g}^{V}
(z^{s};b^r)\; \zeta^{s}\; t^d \;{u^r\over r!}\; \la^{2g-2}\label{last.Hurwitz}
\\
& =& \sum_{\al} N_{d,g}(\al)\;\zeta^\al\; t^d \;{u^r\over r!}\;
       \la^{2g-2}
\nonumber
\end{eqnarray}
where the monomial $\zeta^{\al}$ is dual to $z^{\al}$. (The
$\ell!$ appears because our relative invariant orders the points in
the inverse image of $p$, while the Hurwitz numbers do not.)  This is
a standard generating function for the Hurwitz numbers.

\numbereddemo{Example} The GT invariant of $\P^2 $ relative to a line $L$ is the collection
of enumerative invariants introduced by Caporaso and Harris in
\cite{CH}.
\enddemo

In \cite{CH}, Caporaso and Harris establish a recursion formula for
the number of nodal curves in $\P^2$.  They separate the set of nodal
curves into classes according to how the curves intersect a fixed line
$L$.  Specifically, for each pair of finite sequences $\al=(\al_1,
\al_2,\dots)$ and $\beta=(\beta_1, \beta_2,\dots)$ they consider the
number $N_{d,\delta}(\al,\b)$ of degree $d$ curves with $\delta$
double points, having a contact with $L$ of order $k$ at $\al_k$ fixed
points and at $\b_k$ unspecified points of $L$ for each $k=1,2,\dots$,
and passing through $r=2d+g-1+\ell(\beta)$ fixed points off $L$.  Note
that $\delta$ is determined by the adjunction formula
$2g=(d-1)(d-2)-2\de$.

     From our viewpoint $V=L$ is a symplectic submanifold of $X= \P^2$
with no rim tori.  As in Example \ref{ex9.1} the homology class $A$ of
the map is given by the degree $d$ and we are imposing no constraints
on the complex structure. This time ${\cal S}V$ is the disjoint union
of products of copies of $V=\P^1$.  Since $V$ has only
even-dimensional homology, the relative invariant is again unchanged
under permutations of the marked points that are mapped into $L$.  We
can then associate the generator of $H_0(V_s)$ with the monomial
$y^s=y_{s_1}\cdots y_{s_\ell}$, the generator of $H_2(V_s)$ with
$z^s=z_{s_1}\cdots z_{s_\ell}$. Thus
$$
H_*({\cal S}V)\ =\ \bigoplus_s\ (\Z\oplus \Z)\ = \ \Z[y_1, z_1,y_2, z_2,\dots].
$$
Then
\begin{equation}
N_{d,\delta}(\alpha,\beta)\ =\ {\rm GT}_{\P^2,dL,\chi,r}^{L} (p^r; y^\al, z^\b)\ =\
{\rm GT}_{\P^2,dL,\chi,r}^{L} (p^r;C_{(s,\gamma)})\hskip.4in
\label{last.CH}
\end{equation}
where $\chi=2-2g=-d(d-3)+2\delta$ and where $s$ is any one  of the ordered
sequences such that the basis
element (\ref{last.basis}) satisfies
$$
\al_k\ =\ \mbox{Card}\ \{i\;|\;(s_i,\gamma_i)=(k,[p])\},
$$
and
$$
\b_k\ =\ \mbox{Card}\ \{i\;|\;(s_i,\gamma_i)=(k,[\P^1])\}.
$$

After matching notation through (9.4), we see that the Caporaso-Harris
recursion formula is a consequence of the gluing theorem
for relative invariants; see \cite{IP3}.

\numbereddemo{Example} The GT invariant of the elliptic surfaces $E(n)$ relative to a fiber
$F$.
\enddemo

This example appeared in \cite{IP2}. Here $E(n)\to \P^1$ is the
elliptic surface with a section of self-intersection $-n$, so that  $E(0)=
\P^1\ti T^2$, $E(1)$ is the rational elliptic surface, and $E(2)=K3$,
each regarded as a symplectic manifold.  We focus on counting the
genus 1 (Euler characteristic~$0$) curves representing multiples of the
fiber class.  For generic $(J,\nu)$ that count is given by the numbers
${\rm GT}_{E(n), mF,0}$.  As in \cite{IP2} these agree with the
Seiberg-Witten invariants and are determined by the generating
function
\bear
\sum_m {\rm GT}_{E(n), mF,0}\ t_F^m\ =\ (1-t_F)^{n-2}.
\label{9.GT1}
\eear

The geometric interpretation of this is given in \cite{T}, \cite{IP2}
and \cite{IP3}.  For generic $J$ and $\nu=0$ there are exactly $n-2$
holomorphic fibers; these are of type $(0,-)$ for $n>2$ and type
$(0,+)$ for $n<2$.  The type determines the contribution to the GT
invariant of the maps which multiply cover these fibers when we move
from $(J,0)$ to a generic $(J,\nu)$.  For type $(0,+)$ all covers
contribute, giving the factor $(1-t_F)^{-1}$, while curves of type
$(0,-)$ contribute the opposite factor $(1-t_F)$.

       Now fix a generic fiber $F_0$ and restrict attention to
$F_0$-compatible $(J,0)$.  Then $D^N$ is a complex operator by Lemma
\ref{DNcxlemma}, so that $F_0$ is a holomorphic curve of type $(0,+)$
(cf.\ \cite{T}).  The relative GT invariant then, by the Definition
\ref{defn4.Vregular}, does not contain the contribution of $F_0$ and
its multiple coverings.  Thus
\bear
\sum_m {\rm GT}^F_{E(n), mF,0}\ t_F^m \ =\ (1-t_F)^{n-1}.
\label{9.GT2}
\eear
In particular, the absolute and relative GT invariants are different. In
this case there
are no rim tori, $F_0$-regular maps representing multiples of the fiber
never intersect $F_0$, and (\ref{9.GT2}) agrees with the relative
Seiberg-Witten invariants.

Similarly,  the GT invariant of $E(0)=S^2\ti T^2$ relative to two copies
of $F$ is
$$
\sum_m {\rm GT}^{F,F}_{E(0), mF,0}\ t_F^m\ =\ 1.
$$

 \vglue16pt\centerline{\bf Appendix}

\vglue12pt

The space ${\cal J}^V$ of almost complex structures compatible with
$V$ was defined in Section 3.  Here we show that ${\cal J}^V$ is
nonempty and path-connected.  This fact was used in Section 7 to show
that the relative GW invariants depend only on the symplectic
structure.

An almost
complex structure $J$ on a symplectic manifold $(X,\w)$ is {\it compatible}
with $\w$ if
\bear \speqnu{A.1}
g(X,Y)=\w(X,JY)
\label{A.1}
\eear 
defines a Riemannian metric; this implies that $g(JX,JY)=g(X,Y)$.  Such a
compatible  $J$
can always be constructed, as follows.    After we fix  a ``background''
metric
$g'$, $\w$ defines a skew-symmetric endomorphism
$A$ of
$T^*_p X$ at each point
$p\in X$ by $\w(X,Y)=g'(AX,Y)$.  From linear algebra, any
$A\in\mbox{GL}(n)$ can be {\it uniquely}
expressed as
$A=JS$ where
$J$ is orthogonal  and $S$ is positive definite and symmetric.   Then
$(-J^2)(J^tSJ)=-J(J J^t)SJ=-AJ=A^tJ=S J^tJ=  S$.  Since $J^tSJ$ is
positive definite and symmetric,
the uniqueness of the decomposition gives $J^2=-\mbox{Id}$.  Thus $J$ is an
almost-complex structure, and
then $g(X,Y)=\w(X,JY)$ is a $J$-compatible metric.

Given a symplectic submanifold $V\subset X$, let   $N^g\subset TX$  and
$N^\w\subset TX$ denote the normal bundles to $V$ defined by the metric $g$
and the symplectic
form $\w$ respectively.

\specialnumber{A.1}
\proclaim{Lemma}
\label{N=NLemma}
For compatible $(\w, g, J)${\rm ,} $V$ is $J$\/{\rm -}\/invariant if and only if $N^\w=N^g$.
\endproclaim

\pf  If $N^\w=N^g$ then for any $X\in N^\w$ and $v\in TV$, we have
$g(X,Jv)=-\w(X,v)=0$, so that
$Jv\in TV$; thus $V$ is $J$-invariant.  Conversely,  if $V$ is
$J$-invariant, the equation $g(X,v)=\w(X,Jv)$ implies that $N^\w=N^g$.
\enddemo

\specialnumber{A.2}
\proclaim{Theorem}
\label{JConnectedtheorem}
The space ${\cal J}^V$ of pairs $(J,\nu)$ satisfying {\rm (\ref{defJVeq})}
is nonempty and
path\/{\rm -}\/connected.
\endproclaim

\pf   Replacing $\nu$ by $t\nu$, $0\leq t \leq 1$, gives a retraction of
${\cal J}^V$  to the
space ${\cal J}^V_0$ of $J$  satisfying (\ref{defJVeq}a)   and
(\ref{defJVeq}b).  It therefore
suffices to show that  ${\cal J}^V_0$  is nonempty and  path-connected.

Given triples $(\w, g_0, J_0)$ and $(\w, g_1, J_1)$ in  ${\cal J}^V_0$, we
can apply the above
construction to  the path $h_t=(1-t)g_0+tg_1$ to get a  homotopy  $(\w,
g_t, J_t)$ in which
the decomposition $TV\oplus N^\w$ is orthogonal under  $g_t$ and preserved
by  $J_t$.  Thus
each $J_t$ satisfies   (\ref{defJVeq}a).

To finish the proof, we will modify the path $(\w, g_t, J_t)$ to a path
$(\w, \tilde{g}_t, \tilde J_t)$ which also satisfies
(\ref{defJVeq}b).  For
notational
simplicity we will omit the subscript~$t$.

The Nijenhuis tensor, multiplied by  $-J$, defines a linear map 
$L:N_V\to\mbox{Hom}(TV,
N_V)$  by
\best
L_\xi(v) & = &  \left(\,J[v,\xi] - [v,J\xi] - [Jv,\xi]-J[Jv,J\xi] \,\right)^N\\[6pt]
& = & \left[\, (\nabla_{\xi} J)(v) + J(\nabla_{J\xi}
J)(v) - (\nabla_{v} J)(\xi)  - J(\nabla_{Jv} J)(\xi) \,\right]^N 
\eest
for $\xi\in N_V$.   This is tensorial and $J$-anti-linear in $v$ and 
$\xi$,  and depends
only on $J$ (the second formula above holds for the Levi-Civita connection
of {\it any} Riemannian metric).  Extend $L$ to a map 
$L:N_V\to\mbox{End}(TX)$ by setting
$L_\xi(\eta)=0$ for $\xi,\eta\in N_V$, and let 
$L^t:N_V\to\mbox{End}(TX)$ be its transpose.
By extending $L$ to a neighborhood of $V$,  integrating for a short 
distance along the  lines normal
to
$V$,  and extending  arbitrarily, we can find a 
$K\in\Gamma(\mbox{End}(TX))$ whose 1-jet along $V$
satisfies
\bear \speqnu{A.2}
   K|_V= 0 \qquad\mbox{and}\qquad \nabla_\xi K  = -\frac{1}{2}\left(
L_\xi+L_\xi^t\right) \qquad \forall \xi\in N_V.
\label{JKKJ}
\eear
Then $KJ=-JK$ and $K$ is self-adjoint with respect to $g$. 
Consequently, $JK$ is self-adjoint,
so that
$$
g'(X,Y)\ =\  g\left(e^{KJ}X, Y\right)
$$
defines a Riemannian metric, and it is straightforward to check that 
$J':=e^{JK}J=Je^{KJ}$ is
orthogonal with respect to $g'$.  With $g'$ as  background metric, 
the  procedure described after
(\ref{A.1}) yields a  compatible triple $(\w,\tilde g, \tilde J)$ 
where   $A'=\tilde{J} \tilde{S}$
satisfies
$$
g(JX,Y) \ =\   \w(X,Y)\ =\ g'(A'X,Y) \ =\
g\left( e^{KJ}A'X,Y\right)
$$
and therefore $A'=J'$. The uniqueness of the
factorization $A'=\tilde{J}\tilde{S}=J'\cdot I$ then implies that
$\tilde{J}= J'=J+K+\cdots$ where the dots denote
terms that vanish to second order along $V$.  With that, we can evaluate
$$
\tilde L_\xi(v)\ =\ \left[(\nabla_{\xi} \tilde J)(v) + \tilde
J(\nabla_{\tilde J\xi}
\tilde J)(v) -  (\nabla_{v} \tilde J)(\xi)  - \tilde J
(\nabla_{\tilde Jv} \tilde
J)(\xi)\right]^N
$$
   along $V$.  Using  equation (\ref{JKKJ}), and
the facts that
$L^t_\xi(v)=0$ for \pagebreak $v\in TV$ and $JL_{J\xi}(v)=L_\xi(v)$, we find that
$$
\tilde L_\xi(v)\ =\  L_\xi(v) -\frac12 L_\xi(v) -\frac12 JL_{J\xi}(v)\ =\ 0.
$$
Therefore $(\w,\tilde
g, \tilde J)$ is a  compatible triple satisfying  (\ref{defJVeq}a,b).

Applying this procedure to the path  $(\w, g_t, J_t)$ does not change 
$g_t$ or $J_t$ at
$t=0,1$ (where $L_\xi(v)$ already vanishes) and hence  gives
the desired path $(\w, \tilde{g}_t, \tilde J_t)$.
\enddemo

\AuthorRefNames [IPI]
 
 \end{document}